\documentclass[12pt,reqno]{amsart}
\usepackage{fixltx2e}                    %% Get the latest LaTeX enhancements
\usepackage{mathpazo}  % This provides the palatino font and some nice "old style" numbers.
\usepackage[utf8]{inputenc}            %% Customize to your encoding !!!
\usepackage{amsmath}                     %% Imprescindible
\usepackage{amssymb, latexsym, stmaryrd, amsthm, dsfont, amsfonts, amsbsy,amsthm, amsmath, mathrsfs}            %% Imprescindible
\usepackage{mathtools}                   %% for some math typesetting tricks
\usepackage{bm}                          %% for good bold math symbols
\usepackage{enumerate}                   %% enhanced "enumerate" environments
\usepackage{verbatim}                    %% for enhanced "comment" environment
\usepackage{url}   
\usepackage{lscape}                      %% for \url command in bibliography
\usepackage{nicefrac}   % for nice fractions 1/2, 3/4 etc.

\usepackage[margin=1.2in]{geometry}

\usepackage{scalerel}

\usepackage{microtype}  
\usepackage[all, knot]{xy}
        \xyoption{arc} 
        \xyoption{web}                 %% for the best typographic result
\makeatletter                            %% A correction needed for microtype:
\def\MT@register@subst@font{\MT@exp@one@n\MT@in@clist\font@name\MT@font@list
 \ifMT@inlist@\else\xdef\MT@font@list{\MT@font@list\font@name,}\fi}
\makeatother

\usepackage[pdftex,bookmarks,bookmarksnumbered,linktocpage,   %%  customize to your liking
         colorlinks,linkcolor=blue,citecolor=blue]{hyperref}

\newcommand{\bit}{\begin{itemize}}    % but see also \benbullet below
\newcommand{\eit}{\end{itemize}}
\newcommand{\ben}{\begin{enumerate}}
\newcommand{\een}{\end{enumerate}}
\newcommand{\benroman}{\ben[\normalfont (i)]}  % *
\let\eroman\een
\newcommand{\bde}{\begin{description}}
\newcommand{\ede}{\end{description}}

\newcommand{\mathrmL}{{\mathchoice{\mbox{\rm\L}}{\mbox{\rm\L}}{\mbox{\rm\scriptsize\L}}{\mbox{\rm\tiny\L}}}}

\newcommand{\Rat}{{\mathbb{Q}}}
\newcommand{\Real}{{\mathbb{R}}}

\newcommand{\Prime}{\ensuremath{\sf{Prime}}} 
\newcommand{\?}{\ensuremath{\mkern0.4\thinmuskip}}   % very small math space
\let\models=\vDash                          % models

\let\leq=\leqslant
\let\geq=\geqslant
\let\epsilon=\varepsilon
\let\Lambda\varLambda
\let\Gamma\varGamma
\let\Delta\varDelta
\let\Lambda\varLambda
\let\Omega\varOmega
\let\Theta\varTheta
\let\Xi\varXi
\let\Pi\varPi
\let\Sigma\varSigma
\newcommand\BoldLambda{\ensuremath{\bm{\Lambda}}}

\newcommand{\alg}[1]{ {\ensuremath{\bm {#1}}}}

\let\class=\mathsf                              %  classes (of any things)
\newcommand{\Alg}[1]{\ensuremath{\class {#1}}}	% classes of algebras
\let\oper=\mathbb                               %  operators
\bmdefine{\A}{A}                                %  particular algebras
\bmdefine{\2}{2}
\bmdefine{\B}{B}
\bmdefine{\D}{D}
\bmdefine{\M}{M}                                %  the monoid of substitutions
\bmdefine{\LLL}{L}                              %  algebraic language
\bmdefine{\Fm}{Fm}                              %  formula algebra
\bmdefine{\zerou}{[0{,}1]}  
\bmdefine{\T}{T}                                %  particular algebras

\newcommand\standardL{\ensuremath{{\alg{R}_{\?\text{\sl \L}}^{-}}}}
\newcommand\standardG{\ensuremath{{\alg{R}_{G}^{-}}}}
\newcommand\standardP{\ensuremath{{\alg{R}_{P}^{-}}}}

\newcommand\standardLQ{\ensuremath{{\alg{R}_{\?\text{\sl \L}}}}}
\newcommand\standardGQ{\ensuremath{{\alg{R}_{G}}}}
\newcommand\standardPQ{\ensuremath{{\alg{R}_{P}}}}
\newcommand\standardPQpower[1]{\ensuremath{{\alg{R}_{P}^{#1}}}}

\newcommand\rationalL{\ensuremath{{\alg{Q}_{\?\text{\sl \L}}^{-}}}}
\newcommand\rationalG{\ensuremath{{\alg{Q}_{G}^{-}}}}
\newcommand\rationalP{\ensuremath{{\alg{Q}_{P}^{-}}}}

\newcommand\rationalLQ{\ensuremath{{\alg{Q}_{\?\text{\sl \L}}}}}
\newcommand\rationalGQ{\ensuremath{{\alg{Q}_{G}}}}
\newcommand\rationalPQ{\ensuremath{{\alg{Q}_{P}}}}

\newcommand\prlog{\ensuremath{\mathbf{P}}}
\newcommand\lulog{\ensuremath{\textbf{\L}}}
\newcommand\golog{\ensuremath{\mathbf{G}}}
\newcommand\rprlog{\ensuremath{\mathbf{RP}}}
\newcommand\rlulog{\ensuremath{\textbf{R\L}}}
\newcommand\rgolog{\ensuremath{\mathbf{RG}}}

\newcommand\bllog{\ensuremath{\mathbf{BL}}}

\newcommand{\boland}{\mathop{\wedge\kern-6.4pt\wedge}}
\newcommand{\foo}{\mathop{\bigwedge\kern-11.5pt\bigwedge}} % for display
\newcommand{\boor}{\mathop{\vee\mbox{}\kern-6.4pt\vee}}
\newcommand{\doo}{\mathop{\bigvee\kern-11.5pt\bigvee}}
\newcommand{\boto}{\Longrightarrow}
\newcommand{\boneg}{\text{\raisebox{-.5pt}{$\neg$}\mbox{}\kern-8.7pt\raisebox{1pt}{$\neg$}}}

\newcommand{\VVV}{\oper{V}}                     %  class operators
\newcommand{\QQQ}{\oper{Q}}

\newcommand{\HHH}{\oper{H}}

\newcommand{\PPP}{\oper{P}}

\newcommand{\PPU}{\oper{P}_{\!\textsc{u}}^{}}

\newcommand{\SSS}{\oper{S}}

\newcommand{\III}{\oper{I}}

\newcommand{\Con}{\mathrm{Con}\?}                            % congruences
\newcommand{\Var}{\mathnormal{V\mkern-.8\thinmuskip ar}} % sentential variables
\bmdefine{\boldstar}{\mathchoice{\textstyle*}{\textstyle*}{\textstyle*}{\scriptstyle*}}

\bmdefine{\btau}{\tau}                                  %  transformer tau
\bmdefine{\brho}{\rho}                                  %  transformer rho

\newcommand{\semeq}{\mathrel{=\joinrel\mathrel\vert\mkern0mu\mathrel\vert\joinrel=}}  % semantical equivalence
\bmdefine{\leibniz}{\Omega}        %  Leibniz operator

\bmdefine{\frege}{\Lambda}         %  Frege operator

\makeatletter
\newcommand{\tarskidsp}{\mathord%
   {\m@th\raisebox{0pt}[0pt][0pt]{$\stackrel%
   {\raisebox{-2.7pt}[0ex][0pt]{$\displaystyle \,\?\thicksim$}}%
   {\displaystyle\leibniz}$}}}
\newcommand{\tarskitxt}{\mathord%
   {\m@th\raisebox{0pt}[0pt][0pt]{$\stackrel%
   {\raisebox{-2.7pt}[0ex][0pt]{$\,\?\thicksim$}}{\displaystyle\leibniz}$}}}
\newcommand{\tarskiscr}{\mathord%
   {{\m@th\raisebox{0pt}[0pt][0pt]{$\stackrel%
   {\raisebox{-2.4pt}[0ex][0pt]{$\scriptstyle \,\?\thicksim$}}%
   {\scriptstyle\leibniz}$}}}}
\newcommand{\tarskiscrscr}{\mathord%
   {{\m@th\raisebox{0pt}[0pt][0pt]{$\stackrel%
   {\raisebox{-2pt}[0ex][0pt]{$\scriptscriptstyle \,\?\thicksim$}}%
   {\scriptscriptstyle\leibniz}$}}}}
\newcommand{\tarski}{\@ifnextchar ^ %
   {\mathchoice{\tarskidsp\kern-.07em}{\tarskitxt\kern-.07em}%
   {\tarskiscr\kern-.07em}{\tarskiscrscr\kern-.07em}}%
   {\mathchoice{\tarskidsp}{\tarskitxt}{\tarskiscr}{\tarskiscrscr}}}
\makeatother

\theoremstyle{theorem}
\newtheorem{Theorem}{Theorem}[section]
\newtheorem{Lemma}[Theorem]{Lemma}
\newtheorem{Corollary}[Theorem]{Corollary}
\newtheorem{Proposition}[Theorem]{Proposition}

\newtheorem{Claim}[Theorem]{Claim}

\theoremstyle{definition}
\newtheorem{law}[Theorem]{Definition}

\theoremstyle{remark}

\newtheorem{Remark}[Theorem]{Remark}

\newcommand{\C}{\boldsymbol{C}} 

\usepackage{xcolor}

\begin{document}

\title[Structural completeness in many-valued logics with rational constants]{Structural completeness in many-valued logics with rational constants}

\author{J. Gispert}
\email{jgispertb@ub.edu}
\address{Department of Mathematics and Computer Science, University of Barcelona, Barcelona, Spain}

\author{Z. Hanikov\'a}
\email{hanikova@cs.cas.cz}
\address{Institute of Computer Science of the Czech Academy of Sciences, Prague, Czech Republic}

\author{T. Moraschini}
\email{tommaso.moraschini@ub.edu}
\address{Department of Philosophy, University of Barcelona, Barcelona, Spain}

\author{M. Stronkowski}
\email{m.stronkowski@mini.pw.edu.pl}
\address{Institute of Computer Science of the Czech Academy of Sciences, Prague, Czech Republic}
\address{Faculty of Mathematics and Information Science, Warsaw University of Technology}

\date{\today}
\keywords{Structural completeness, admissible rule, quasivariety, fuzzy logic, product logic, {\L}ukasiewicz logic, G\"odel logic, rational Pavelka logic}
\renewcommand{\thefootnote}{\fnsymbol{footnote}} 
\renewcommand{\thefootnote}{\arabic{footnote}} 
\maketitle

% 2010 Math subject classification
% 03B52 Fuzzy logic
% 06D35 MV-algebras
% 06D20 Heyting algebras
% 08C15 Quasivarieties

%\tableofcontents

\begin{abstract}
The logics $\rlulog$, $\rprlog$, and $\rgolog$ have been obtained by expanding
{\L}ukasiewicz logic $\lulog$, product logic $\prlog$,  and G\"odel--Dummett logic $\golog$ 
with rational constants.
We study the lattices of extensions and structural completeness of
these three expansions, obtaining results that stand in
contrast to the known situation in $\lulog$, $\prlog$, and $\golog$. 
Namely, $\rlulog$ is hereditarily structurally complete.
$\rprlog$ is algebraized by the variety of rational product algebras that we
show to be $\mathcal{Q}$-universal. We provide a base of admissible rules in $\rprlog$,
 show their decidability, and characterize passive structural completeness
for extensions of $\rprlog$.  Furthermore, structural completeness, hereditary structural
completeness, and active structural completeness coincide for extensions of $\rprlog$, and
this is also the case for extensions of $\rgolog$, where in turn passive
structural completeness is characterized by the equivalent algebraic semantics
having the joint embedding property. For nontrivial axiomatic extensions of
$\rgolog$  we provide a base of admissible rules. We leave the problem open
whether the variety of rational G\"odel algebras is $\mathcal{Q}$-universal.

\end{abstract}

\section{Introduction}

This work brings together two lines of research:
admissible rules and lattices of extensions of logics on the one side, and 
propositional fuzzy logic with constants for rational numbers on the other. Either of these lines is
native to nonclassical
logics and trivializes in the classical case. 

In the realm of admissibility, it is common to identify \emph{logics} with finitary substitution invariant consequence relations $\vdash$ on the set of formulas of some algebraic language. Formulas $\varphi$ such that $\emptyset \vdash \varphi$ are then called \emph{theorems} of $\vdash$. A rule $\gamma_1, \dots, \gamma_n \rhd \varphi$ is \textit{derivable} in a logic $\vdash$ when $\gamma_1, \dots, \gamma_n \vdash \varphi$. It is \emph{admissible} in $\vdash$ provided that the set of theorems of $\vdash$ 
is closed under that rule.  The derivability of a rule entails its admissibility in $\vdash$, but the converse fails in general.  Indeed, the \textit{structural completion} of a logic $\vdash$ is the only logic whose derivable rules are precisely the rules that are admissible in $\vdash$. Because of this, questions typically asked about
derivability, such as finding an axiomatization or settling decidability, pertain
also to admissibility. 
In a \emph{structurally complete} logic, the set of 
admissible rules coincides with the set of derivable rules,
and a logic is called \emph{hereditarily} structurally complete 
if the property of structural completeness is shared by all of its extensions.

Research in structural completeness is well established in intermediate and transitive modal logics.\
% apparently started with Lorenzen in 1955(!) what's in that paper? 
Landmarks include the work of Rybakov on decidability of admissible rules, 
covered in the monograph \cite{Ry97}, 
% the problem of decidability stated by Friedman 1975
Ghilardi's investigation of the relationship  of admissibility to unification \cite{ghilardi:unificationInt}, 
and Iemhoff's construction of explicit bases for the rules admissible in the intuitionistic logic \cite{Ie01},
independently discovered by Rozi\`{e}re \cite{Roziere93}.  Hereditarily structurally complete logics have been described in the realm of intermediate logics by Citkin \cite{Citkin78a,Citkin87}, see also \cite{BezMor19,Citkin19HSC}, and by Rybakov in that of transitive modal logics \cite{Ryb95} (see also \cite{Dziobiak}). Structural completeness is language sensitive:
the pure implication fragment of the intuitionisitc logic has been known to be hereditarily structurally complete due to Prucnal's work \cite{Prucnal:StructuralCompleteness}, 
yet the implication-negation fragment is incomplete \cite{CintulaMetcalfe:AdmissibleImpleNegInt}.

During the last two decades, research in structural completeness has turned also to the family of fuzzy logics.
Three logics of interest in this paper --- {\L}ukasiewicz logic $\lulog$, product logic $\prlog$, 
and G\"odel--Dummett logic $\golog$ --- can be obtained as axiomatic extensions of H\'ajek's basic logic $\bllog$ \cite{Ha98},
even if they were defined independently prior to the definition of $\bllog$.
{\L}ukasiewicz logic was first introduced  in \cite{Lukasiewicz-Tarski:Untersuchungen}.
Finite-valued semantics for G\"odel--Dummett logic\footnote{Henceforth we write just \textit{G\"odel
logic}, as is common in the referenced literature.} was considered in G\"odel's analysis
of intuitionistic logic \cite{Godel:ZumAussagen}, 
while Dummett provided an axiomatization of the infinite-valued case in \cite{Dummett:GodelLogic}.\ Product logic first appeared in \cite{Hajek-Godo-Esteva:ProductLogic}; see also \cite{Ha98}.\ All logics in the $\bllog$ family are algebraizable in the sense of Blok and Pigozzi \cite{BP89}:
% their equivalent semantics is the class of BL-algebras, i.e.,  bounded, commutative, integral, semilinear, and % divisible Full Lambek (FL) algebras as presented in \cite{Galatos-JKO:ResiduatedLattices}. 
the equivalent algebraic semantics of the three logics are the varieties of MV-algebras, product algebras, and G\"odel algebras respectively.
MV-algebras and product algebras each have a specific tight connection to lattice-ordered abelian groups, while
G\"odel algebras coincide with Heyting algebras in which the equation $(x\to y)\lor(y\to x) \thickapprox 1$ holds.

While $\golog$ and $\prlog$ are hereditarily structurally complete \cite{DzWr73,CiMe09},  $\lulog$  is structurally incomplete \cite{Dz08} and a base for its admissible rules was exhibited by Je\v r\'abek  \cite{Jerabek:BasesAdmissibleRulesLuk}, see also  \cite{Jerabek:AdmissibleRulesLuk,Jerabek:ComplexityAdmissibleRulesLuk}. Admissibility in extensions of $\lulog$ was investigated in \cite{Gisp16,Gisp17}.\footnote{More generally, results
on logics without weakening have been obtained in \cite{OlRaAl08,RaftSwiry16} . } 
 Finally, works addressing variants of structural completeness, such as \emph{active}
and \emph{passive} rules also studied in this paper, 
include  \cite{DzSt201x,Gisp18,MeRo13a,MorRafWan19PSC,RaxxNJ,Wro09}.  

The aim of this paper is to take this line of research further and  
look at structural completeness for expansions of {\L}ukasiewicz, product,  and G\"odel 
logic with rational constants. 
The pedigree of these logics goes back to the pioneering works of Goguen \cite{Go69} and Pavelka \cite{Pa79a,Pa79b,Pa79c}.
Expanding the language with constants can be viewed as taking advantage of 
the rich algebraic setting to gain more expressivity; see, e.g., \cite{Belohlavek:PavelkaRetrospectProspect,Cintula:noteAxiomatizationsPavelka,Esteva-Godo-Montagna:LPi,Esteva-Godo-Noguera:ExpandingWNM,Hanikova:ComplexityValidDegrees,TaTi92}.

More specifically, while the logics live in an ambience of many truth values, 
they formally derive only statements that  are fully true. 
It has been established by Goguen and Pavelka that if constants with suitable axioms  
are added to a logic such as $\lulog$, 
one can walk around this limitation by employing formulas of the form $\boldsymbol{c}\to\varphi$,
with $\boldsymbol{c}$ a constant and $\varphi$ any formula of the language:  
assignments sending $\boldsymbol{c}\to\varphi$ to the top element are precisely those that 
send $\varphi$ in the upset of the value of $\boldsymbol{c}$.   
Via this simple hedging device, the existing deductive machinery of the logic 
(still ostensibly focused on fully true statements) 
 enables deduction on graded statements.

The version of $\lulog$ with rational constants in \cite{Hajek:ArHierarchyI,Ha98}
 has become known as \emph{rational Pavelka logic}. 
Here we refer to this logic as \textit{rational {\L}ukasiewicz logic} ($\rlulog$) 
to have uniform names over all three expansions, the other two being
\textit{rational product logic} $\rprlog$ and \textit{rational G\"odel logic} $\rgolog$ \cite{EsGiGoNo07,Esteva-Godo-Noguera:Rational,Hajek:ComplexityRational,Savicky-CEGN:ProductTruthConstants}. These logics are algebraized, respectively, by the varieties of \textit{rational MV-algebras}, \textit{rational product algebras}, and \textit{rational G\"odel algebras}.

For each of the three logics, we provide information on the lattice of its extensions and 
identify the structurally complete ones.\
The following striking reversal in the lattice structure of extensions instantiates
the already mentioned language sensitivity of the notions studied in this paper.  
While $\lulog$ is known to be structurally incomplete and the lattice of its extensions is 
dually isomorphic to lattice of quasivarieties of the $\mathcal{Q}$-universal variety of MV-algebras \cite{AdaDzi94}\footnote{Even if the class of MV-algebras is $\mathcal{Q}$-universal, insights into the structure of quasivarieties generated by MV-chains were 
provided in \cite{Gp02,GpMuT99,GpT98}.}, 
$\rlulog$ is hereditarily structurally complete, there being no consistent extensions.
On the other hand, the lattice of extensions of $\rprlog$ is dually isomorphic to the lattice 
of subquasivarieties of rational product algebras, which we show to be $\mathcal{Q}$-universal,
and the only structurally complete extensions are the logic of the rational product algebra on the rationals in $[0,1]$ with the natural order and the three proper axiomatic extensions of $\rprlog$ term-equivalent to the extensions of $\prlog$. This contrasts with the known situation in $\prlog$,
which is hereditarily structurally complete and whose lattice of extensions is a three-element chain.\ Lastly, while $\golog$ is hereditarily structurally complete and has denumerably many extensions, $\rgolog$ is structurally incomplete and has a continuum of axiomatic extensions.

The paper is structured as follows. Sections \ref{Sec:2}, \ref{Sec:3}, and \ref{Sec:4} review the rudiments of the theory of 
quasivarieties, structural completeness, and fuzzy logics respectively. 
Section \ref{Sec:5} establishes the $\mathcal{Q}$-universality of the class of rational product algebras. 
Section \ref{Sec:6} is dedicated to structural completeness results in extensions of $\rprlog$;
in particular, Theorem  \ref{Thm:base-product} provides a base of rules admissible  in $\rprlog$,
while Corollary \ref{Cor:decid-RPLOG-admiss} establishes their decidability. Theorem \ref{Thm:SC-product} characterizes
structurally complete extensions of $\rprlog$; it turns out that structural completeness,
active structural completeness, and hereditary structural completeness coincide for extensions of $\rprlog$
(Corollaries \ref{Cor:HSC-product} and \ref{Cor:ASC-product}).
Corollary \ref{Cor:PSC-product} offers a characterization of passively structurally complete extensions. 
Section \ref{Sec:7} studies the lattice of extensions of $\rgolog$:
by Corollary \ref {Cor:axiomatic-extensions-RG}, 
already the lattice of axiomatic extensions of $\rgolog$ is an uncountable chain.\
The lattice of $\rgolog$-extensions is also easily seen to have uncountable antichains; but 
we do not know whether the class of rational G\"odel algebras might be $\mathcal{Q}$-universal.
Section \ref{Sec:8} studies structural completeness in $\rgolog$; for extensions of this logic, 
structural completeness, hereditary structural completeness, and active structural completeness coincide
with the extension being algebraized by a quasivariety generated by a chain in $\class{RGA}$
 (Theorem \ref{Thm:RGA-SC}); moreover, all such rational G\"odel chains are
characterized. Passive structural completeness for $\rgolog$-extensions is characterized in terms of
being algebraized by a quasivariety having the joint embedding property (Theorem \ref{Thm:RGA-PSC}).
Theorem \ref{Thm:SCbase-godel} provides a base of admissible rules in any nontrivial axiomatic extension 
of $\rgolog$.
% decidability?
Finally,  Theorem \ref{thm_structural_completeness_RMV} in section \ref{Sec:9} shows that $\rlulog$  lacks proper consistent extensions and, therefore, is hereditarily structurally complete.

\section{Varieties and quasivarieties}\label{Sec:2}

% \marginpar{change (quasi)varieties to (quasi)equational classes and 
% define (quasi)varieties as closed under class operators.}

A \textit{quasivariety} is a class of algebras that can be axiomatized by \text{quasiequations}, 
i.e., sentences of the form
\[
\forall \vec{x}((\varphi_{1} \thickapprox \psi_{1} \boland \dots \boland \varphi_{n} \thickapprox \psi_{n}) \boto \varphi \thickapprox \psi).
\]
We admit the case where the antecedent of the above implication is empty, whence universally quantified equations are special cases of quasiequations. Similarly, a \textit{variety} is a class of algebras that can be axiomatized by universally quantified equations, while a \textit{universal class}
is one that can be axiomatized by universally quantified open formulas. 
% While displaying quasiequations, 
It is common to drop the universal quantifiers in the prefix
and work with open formulas.

For a general introduction to the theory of these classes we refer the reader to \cite{Be11g,BuSa00,Go98a} and, in what follows, we shall review some fundamental material only. Varieties, quasivarieties and universal classes 
can be characterized in terms of model-theoretic constructions. Let $\III, \HHH, \SSS, \PPP$, and $\PPU$ be the class operators of closure under isomorphism, homomorphic images, subalgebras, direct products,  and ultraproducts respectively. We assume direct products and ultraproducts of empty families of algebras are trivial algebras. Then, a class of similar algebras $\class{K}$ is a variety precisely when it is closed under $\HHH, \SSS$, and $\PPP$  \cite[Thm.\ I.11.9]{BuSa00}, it is a quasivariety precisely when it is closed under $\III, \SSS, \PPP$, and $\PPU$ \cite[Thm.\ V.2.25]{BuSa00} (see also \cite[Cor.\ 2.4]{Stronkowski:Axiomatizations}),
and it is a universal class precisely when it is closed under $\III$, $\SSS$\, and $\PPU$ \cite[Thm.\ V.2.20]{BuSa00}. Given a class of similar algebras $\class{K}$, the smallest variety and quasivariety 
containing $\class{K}$ will be denoted by $\VVV(\class{K})$ and $\QQQ(\class{K})$, respectively. 
It turns out that
$\VVV(\class{K}) = \HHH\SSS\PPP(\class{K})$ and  $\QQQ(\class{K}) = \III\SSS\PPP\PPU(\class{K})$.
Moreover, the smallest universal class containing $\class{K}$ is $\III\SSS\PPU(\class{K})$.

A \textit{finite partial subalgebra} $\C$ of an algebra $\A$ is a finite subset $C$ of $A$ endowed with the restriction of finitely many basic operations of $\A$.  Given two similar algebras $\A$ and $\B$, a finite partial subalgebra $\C$ of $\A$ is said to \textit{embed} into $\B$ if there exists an injective map $h: C \to B$  such that for every basic $n$-ary partial operation $f$ of $\C$ and $c_1,\dots,c_n \in C$ such that $f^{\A}(c_1,\dots, c_n)\in C$,  we have
\[
h (f^{\A} (c_1,\dots,c_n)) = f^{\B}(h(c_1),\dots,h(c_n)).
\]  
In this case, we say that $h$ is an \textit{embedding} of $\C$ into $\B$.\ When every finite partial subalgebra of $\A$ embeds into $\B$, we say that $\A$ \textit{partially embeds} into ${\B}$. 
% Recall that $\A$ \emph{locally embeds} into  $\B$ if for every finite subset $C$ of $A$ 
% and a finite set $F$ of operation symbols in the language, 
% there exists an injective mapping $e\colon C\to B$ such that for every 
% $n$-ary $f\in F$, $a_1,\ldots,a_n,a\in C$ we have $f^\B(e(a_1),\ldots,e(a_n))=e(a)$ 
% whenever $f^\A(a_1,\ldots,a_n)=a$. 
Partial embeddability is strictly 
connected with universal classes,  because an algebra $\A$ partially embeds into an algebra $\B$ if an only if $\A$ validates the universal theory of $\B$.
% embeddings given in \cite{Ma73d} is slightly different although equivalent to the one given here.

Consider a quasivariety $\class{K}$ and an algebra $\A \in \class{K}$. A congruence $\theta$ of $\A$ is said to be a $\class{K}$-\textit{congruence} if $\A / \theta \in \class{K}$. When ordered under inclusion, the set of $\class{K}$-congruences of $\A$ is an algebraic lattice, which we denote by $\Con_{\class{K}}\A$. On the other hand, the lattice of all congruences of $\A$ will be denoted by $\Con \A$. A congruence of $\A$ is said to be \textit{nontrivial} if it differs form the total relation $A \times A$ and the identity relation $\textup{Id}_{A}$. %Congruences classes will be called \textit{blocks}, 
The kernel of a homomorphism $f$ will be denoted by $\textup{Ker}(f)$. Given $a, c  \in A$, the $\class{K}$-congruence of $\A$ generated by $\langle a, c\rangle$ is denoted by $\textup{Cg}_{\class{K}}^{\A}(a, c)$.

Given a quasivariety $\class{K}$, an algebra $\A \in \class{K}$ is said to be \textit{relatively subdirectly irreducible} (resp.\ \textit{relatively finitely subdirectly irreducible}) in $\class{K}$ if $\textup{Id}_{A}$ is completely meet-irreducible (resp.\ meet-irreducible) in $\Con_{\class{K}}\A$. When $\class{K}$ is a variety, $\Con_{\class{K}}\A = \Con\A$ and $\A$ is said to be simply \textit{subdirectly irreducible} (resp.\ \textit{finitely subdirectly irreducible}).\ The class of algebras that are relatively subdirectly irreducible (resp.\ relatively finitely subdirectly irreducible) in $\class{K}$ will be denoted by $\class{K}_{\textup{RSI}}$ (resp.\ $\class{K}_{\textup{RFSI}}$).  It is well known that every member of a quasivariety $\class{K}$ is isomorphic to a subdirect product of algebras in $\class{K}_{\textup{RSI}}$ \cite[Thm.\ 3.1.1]{Go98a}.  Accordingly, to prove that two quasivarieties $\class{K}$ and $\class{K'}$ are equal, it suffices to show that $\class{K}_{\textup{RSI}} = \class{K'}_{\textup{RSI}}$.

Given a quasivariety $\class{K}$, we denote by $\mathcal{Q}(\class{K})$ the lattice of subquasivarieties of $\class{K}$. On the other hand, a class $\class{V} \subseteq \class{K}$ is said to be a \textit{relative subvariety} of $\class{K}$ if it can be axiomatized by equations relative to $\class{K}$. The lattice of relative subvarieties of $\class{K}$ will be denoted by $\mathcal{V}(\class{K})$. Notice that, when $\class{K}$ is a variety, $\mathcal{V}(\class{K})$ is the lattice of subvarieties of $\class{K}$. A quasivariety $\class{K}$ is said to be \textit{primitive} when all its subquasivarieties are relative subvarieties.

\begin{Theorem}[\protect{\cite[Prop.\ 5.1.22]{Go98a}}]\label{Thm:primitive-distributive}
If $\class{K}$ is a primitive quasivariety, then $\mathcal{Q}(\class{K})$ is a distributive lattice. 
\end{Theorem}

A quasivariety $\class{K}$ has the \textit{joint embedding property} (JEP) when every two nontrivial members $\A$ and $\B$ of $\class{K}$ can be embedded into a common $\C \in \class{K}$. While every variety is generated by its denumerably generated free algebra, it is not true that every variety is generated by a single algebra as a quasivariety. This makes the next result from \cite{Mal66ail} interesting in the context of varieties as well.

\begin{Proposition}[\protect{\cite[Prop.\ 2.1.19]{Go98a}}]\label{Prop:Maltsev-JEP}
A quasivariety has the JEP if and only if it is generated by a single algebra as a quasivariety.
\end{Proposition}

Finally, a quasivariety $\class{K}$ is said to be \textit{Q-universal} if $\mathcal{Q}(\class{M}) \in \HHH\SSS(\mathcal{Q}(\class{K}))$, for every quasivariety $\class{M}$ in a finite language.\footnote{The usual definition of a Q-universal quasivariety $\class{K}$ demands that $\class{K}$ has finite language. In this paper we drop this requirement, because we deal with quasivarieties whose language is always infinite.}  As lattices of quasivarieties in a finite language may be uncountable and need not validate any nontrivial lattice equation \cite{GurTum80,Tumanov88}, the next result follows.

\begin{Proposition}\label{Prop:Q-univeresal-property}
If $\class{K}$ is a $\mathcal{Q}$-universal quasivariety, then $\mathcal{Q}(\class{K})$ is uncountable and does not validate any nontrivial lattice equation.
\end{Proposition}

\section{Structural completeness}\label{Sec:3}

Let $\Var = \{ x_{n} : n \in \omega \}$ be a denumerable set of variables. Given an algebraic language $\mathscr{L}$, we denote by $Fm_{\mathscr{L}}$ the set of formulas of $\mathscr{L}$ with variables in $\Var$. When $\mathscr{L}$ is clear from the context, we shall write $Fm$ instead of $Fm_{\mathscr{L}}$. A (propositional) \textit{logic} $\vdash$ is then a consequence relation on the set of formulas $Fm$ of some algebraic language that, moreover, is \textit{substitution invariant} in the sense that for every substitution $\sigma$ on $Fm$ and every $\Gamma \cup \{ \varphi \} \subseteq Fm$,
\[
\text{if }\Gamma \vdash \varphi \text{, then }\sigma[\Gamma] \vdash \sigma(\varphi).
\]
Furthermore, in this paper logics $\vdash$ are assumed to be \textit{finitary}, i.e., such that
\[
\text{if }\Gamma \vdash \varphi\text{, then }\Delta \vdash \varphi \text{ for some finite }\Delta \subseteq \Gamma. 
\]

Given two logics $\vdash$ and $\vdash'$ such that the language of $\vdash'$ extends that of $\vdash$, we say that $\vdash'$ is an \textit{expansion} of $\vdash$ if, for every set of formulas $\Gamma \cup \{ \varphi \}$ in the language of $\vdash$,
\[
\Gamma \vdash \varphi \Longleftrightarrow \Gamma \vdash' \varphi.
\]
Similarly, given two logics $\vdash$ and $\vdash'$ in the same language, $\vdash'$ is said to be an \textit{extension} of $\vdash$ when $\Gamma \vdash'\varphi$, for every $\Gamma \cup \{ \varphi \} \subseteq Fm$ such that $\Gamma \vdash \varphi$. An extension $\vdash'$ of $\vdash$ is said to be \textit{axiomatic} when there is a set $\Sigma \subseteq Fm$ closed under substitutions such that for all $\Gamma \cup \{ \varphi \} \subseteq Fm$,
\[
\Gamma \vdash' \varphi \Longleftrightarrow \Gamma \cup \Sigma \vdash \varphi.
\]

We shall now review the rudiments of the theory of admissible rules. For a systematic treatment, the reader may consult \cite{Pogorzelski-Wojtylak:Completeness,Ry97}. A formula $\varphi$ is said to be a \textit{theorem} of a logic $\vdash$ if $\emptyset \vdash \varphi$. Moreover, a \textit{rule} is an expression of the form $\Gamma \rhd \varphi$, where $\Gamma \cup \{ \varphi \} \subseteq Fm$ is a finite set. When $\Gamma = \{ \gamma_{1}, \dots, \gamma_{n} \}$, we shall sometimes write $\gamma_{1}, \dots, \gamma_{n} \rhd \varphi$ instead of $\Gamma \rhd \varphi$. A rule $\Gamma \rhd \varphi$ is said to be  \textit{derivable} in a logic $\vdash$ when $\Gamma \vdash \varphi$. It is \textit{admissible} in $\vdash$ when for every substitution $\sigma$ on $Fm$,
\[
\text{if }\emptyset \vdash \sigma(\gamma)\text{ for all }\gamma \in \Gamma, \text{ then }\emptyset \vdash \sigma(\varphi).
\]
In other words, a rule is admissible in $\vdash$ when its addition to $\vdash$ does not produce any new theorem. Clearly, every rule that is derivable in $\vdash$ is also admissible in $\vdash$. If the converse holds, $\vdash$ is said to be \textit{structurally complete} (SC). Logics whose extensions are all structurally complete have been called \textit{hereditarily structurally complete} (HSC). 

Every logic admits a canonical structurally complete extension,  see, e.g., \cite[Lem.\ 1.76 \& Thms.\ 1.78 \& 1.79]{Ry97}.

\begin{Proposition}
Every logic $\vdash$ has a unique structurally complete extension $\vdash^{+}$ with the same theorems. Furthermore, a rule is derivable in $\vdash^{+}$ precisely when it is admissible in $\vdash$.
\end{Proposition}
\noindent In view of the above result, $\vdash^{+}$ has been called the \textit{structural completion} of $\vdash$. Since the derivable rules of $\vdash^{+}$ coincide with those admissible in $\vdash$, a set $\Sigma$ of rules is said to be a \textit{base for the admissible rules on $\vdash$} if its addition to $\vdash$ axiomatizes $\vdash^{+}$.

Structural completeness can be split in two halves.  A rule $\Gamma \rhd \varphi$ is said to be \textit{active} in a logic $\vdash$ if there exists a substitution $\sigma$ such that $\emptyset \vdash \sigma(\gamma)$ for all $\gamma \in \Gamma$. It is said to be \textit{passive} in $\vdash$ otherwise. Then, a logic $\vdash$ is called \textit{actively structurally complete} (ASC) if every active rule that is admissible in $\vdash$ is also derivable in $\vdash$, see \cite{DzSt201x,MeRo13a} (where the adjective \textit{almost} is used instead). Notice that every passive rule is vacuously admissible. Accordingly, $\vdash$ is said to be \textit{passively structurally complete} (PSC) \cite{Wro09} if all rules that are passive in $\vdash$ are also derivable in $\vdash$.

A logic $\vdash$ is \textit{algebraized} by a quasivariety $\class{K}$ \cite{BP89} when there are a finite set of equations $\btau(x)$ and a finite set of formulas $\Delta(x, y)$ such that for every $\Gamma \cup \{ \varphi \} \subseteq Fm$,
\begin{align*}
\Gamma \vdash \varphi &\Longleftrightarrow \bigcup \{ \btau(\gamma) : \gamma \in \Gamma \} \vDash_{\class{K}} \btau(\varphi)\\
x \thickapprox y &\semeq_{\class{K}} \bigcup \{ \btau(\delta) : \delta \in \Delta(x, y) \}
\end{align*}
where $\vDash_{\class{K}}$ is the equational consequence relative to $\class{K}$ \cite{BP89,AAL-AIT-f}.  In this case, $\class{K}$ is uniquely determined \cite[Thm.\ 2.15]{BP89} and is called the \textit{equivalent algebraic semantics} of $\vdash$.

When a logic $\vdash$ is algebraized by a quasivariety $\class{K}$, structural completeness and its variants admit the following purely algebraic characterization, in which $\Fm_{\class{K}}(\omega)$ and $\Fm_{\class{K}}(0)$ denote, respectively, the denumerably and zero-generated free algebras of $\class{K}$.

\begin{Theorem}\label{Thm:Bergman} If a logic $\vdash$ is algebraized by a quasivariety {$\class {K}$}, then

\benroman
\item\label{item:Ber-1} $\vdash$ is SC if and only if $\class{K}$ is generated as a quasivariety by $\Fm_{\class{K}}(\omega)$;
\item\label{item:Ber-2} $\vdash$ is HSC if and only if $\class{K}$ is primitive;
\item\label{item:Ber-3} $\vdash$ is  PSC if and only if every positive existential sentence is either true in all nontrivial members of $\class{K}$ or false in all of them;
\item\label{item:Ber-4} $\vdash$ is ASC if and only if $\A \times \Fm_{\class{K}}(\omega) \in \QQQ (\Fm_{\class{K}}(\omega))$ for every relatively subdirectly irreducible algebra $\A \in \class{K}$. If there is a constant symbol in the language, then we can replace ``$\A \times \Fm_{\class{K}}(\omega) \in \QQQ (\Fm_{\class{K}}(\omega))$'' by ``$\A \times \Fm_{\class{K}}(0) \in \QQQ (\Fm_{\class{K}}(\omega))$'' in this statement.
\eroman
\end{Theorem}
In the above result, items (\ref{item:Ber-1}) and (\ref{item:Ber-2}) are essentially \cite[Props.\ 2.3 \& 2.4(2)]{Be91b}, while (\ref{item:Ber-3}) is \cite[Cor.\ 3.2]{DzSt201x}. Lastly, (\ref{item:Ber-4}) was essentially proved in \cite{DzSt201x}, but see also \cite[Thm.\ 7.3]{RaxxNJ}.

When a logic $\vdash$ is algebraized by a quasivariety $\class{K}$ by means of finite sets of equations and formulas $\btau$ and $\Delta$, the lattice of extensions of $\vdash$ is dually isomorphic to $\mathcal{Q}(\class{K})$ \cite[Cor.\ 3.40]{AAL-AIT-f}. The dual isomorphism is given  by the map that sends an extension $\vdash'$ to the quasivariety axiomatized by the quasiequations
\[
\foo \btau(\gamma_{1}) \boland \dots \boland \foo \btau(\gamma_{n}) \boto \epsilon \thickapprox \delta,
\]
where $\gamma_{1}, \dots, \gamma_{n} \vdash' \varphi$ and $\epsilon \thickapprox \delta \in \btau(\varphi)$. The inverse of this dual isomorphism sends a quasivariety $\class{M} \in \mathcal{Q}(\class{K})$ to the logic axiomatized by the rules
\[
\Delta(\varphi_{1}, \psi_{1}) \cup \dots \cup \Delta(\varphi_{n}, \psi_{n}) \rhd \delta
\]
where $\class{M} \vDash (\varphi_{1} \thickapprox \psi_{1} \boland \dots \boland \varphi_{n} \thickapprox \psi_{n}) \boto \varphi \thickapprox \psi$ and $\delta \in \Delta(\varphi, \psi)$. Furthermore, the dual isomorphism restricts to one between the lattice of axiomatic extensions of $\vdash$ and $\mathcal{V}(\class{K})$. Accordingly, the lattice of extensions (resp.\ axiomatic extensions) of $\vdash$ can be studied through the lens of $\mathcal{Q}(\class{K})$ (resp.\ $\mathcal{V}(\class{K})$). The effect of structural completeness on the lattice of extensions of $\vdash$ is captured by the following results, the first of which is a direct consequence of Theorem \ref{Thm:primitive-distributive} and Theorem \ref{Thm:Bergman}(\ref{item:Ber-2}).

\begin{Corollary}\label{Cor:HSC-distributive}
If an HSC logic $\vdash$ is algebraized by a quasivariety $\class{K}$, then the lattice of extensions of $\vdash$ and $\mathcal{Q}(\class{K})$ are distributive.
\end{Corollary}

\begin{Proposition}[\protect{\cite[Thm.\ 4.3 \& Rmk.\ 5.13]{MorRafWan19PSC}}]\label{Prop:PSC-JEP}
Let $\vdash$ be a logic algebraized by a quasivariety $\class{K}$.\ If $\vdash$ is PSC, then every member of $\mathcal{Q}(\class{K})$ has the JEP. Moreover, for every extension $\vdash'$ of $\vdash$ there exists an algebra $\A$ such that, for every $\Gamma \cup \{ \varphi\} \subseteq Fm$,
\[
\Gamma \vdash' \varphi \Longleftrightarrow  \btau[\Gamma] \vDash_{\A} \btau(\varphi),
\]
where $\btau$ is the set of equations witnessing the algebraization of $\vdash$.
\end{Proposition}

\section{Fuzzy logic}\label{Sec:4}
\label{sec:: Fuzzy logic}

A \textit{BL-algebra} is a structure $\A = \langle A; \land, \lor, \cdot, \to, 0, 1 \rangle$ that comprises a bounded lattice $\langle A; \land, \lor, 0, 1 \rangle$ and a commutative monoid $\langle A; \cdot, 1 \rangle$ such that, for every $a, b, c \in A$,  the \textit{residuation law}
\[
a \cdot b \leq c \Longleftrightarrow a \leq b \to c
\]
holds and 
\[
(a \to c) \lor (c \to a ) = 1 \, \, \, \, \text{ and }\, \, \, \,  a \land c = a\cdot (a \to c).
\]

It follows that the lattice reduct of $\A$ is distributive; see Corollary \ref{Cor:prelinear} below. 
Totally ordered algebras are referred to as \emph{chains}.
Furthermore, the lattice operations can be defined in terms of $\cdot$ and $\to$.  For $\land$ this is a consequence of the above display, while for $\lor$ we have
\[
x \lor y \coloneqq ((x \to y) \to y) \land ((y \to x) \to x) ).
\]
From a logical standpoint, the class of BL-algebras forms a variety that algebraizes H\'ajek's \textit{basic logic} $\mathbf{BL}$ \cite{Ha98}.

Given a BL-algebra $\A$, a nonempty set $F \subseteq A$ is said to be a \textit{filter} of $\A$ if it is upward closed, in the sense that if $a \in F$ and $a\leq c$, then $c \in F$ (an upset),  and it is closed under multiplication,  that is,  if $a, c \in F$, then $a \cdot c \in F$.  A filter $F$ of $\A$ is called \textit{prime} when, for every $a, c \in A$,
\[
\text{if }a \lor c \in F \text{, then }a \in F \text{ or }c \in F.
\] 
When ordered under the inclusion relation, the set $\textup{Fi}\?\A$ of filters of $\A$ becomes a lattice that, moreover, is isomorphic to $\Con\A$.

\begin{Theorem}[\protect{\cite[Lem.\ 2.3.14]{Ha98}}]\label{Thm:prelinear}
\label{Thm: filters and congruences}
Let $\A$ be a BL-algebra.  The map $\theta_{(-)} \colon \textup{Fi}\?\A \to \Con\A$, defined by the rule
\[
\theta_F \coloneqq \{ \langle a, c \rangle \in A \times A : a \to c, c \to a \in F \},
\]
is a lattice isomorphism.  Furthermore, the following conditions are equivalent for a filter $F$ of $\A$:
\benroman
\item $F$ is prime;
\item $\A/\theta_F$ is a chain;
\item $\A/\theta_F$ is finitely subdirectly irreducible.
\eroman
\end{Theorem}
\noindent Henceforth, we will write $\A / F$ as a shorthand for $\A / \theta_F$.

The following observation is instrumental to prove the existence of prime filters in BL-algebras.  Its proof is a straightforward adaptation of \cite[Lem.\ 2.3]{Kihara-Ono:variable}.

\begin{Lemma}
\label{lem:: prime filter lemma}
Let $\A$ be a BL-algebra and $I \subseteq A \smallsetminus \{ 1 \}$ such that $a \lor c\in I$, whenever $a, c\in I$. Then there is a prime filter $F$ of $\A$ disjoint from $I$.
\end{Lemma}

In view of the subdirect decomposition theorem \cite[Thm.\ 3.24]{Be11g},  the second part of Theorem \ref{Thm: filters and congruences} implies the following.

\begin{Corollary}\label{Cor:prelinear}
Every BL-algebra is isomorphic to a subdirect product of BL-chains. As a consequence,  the lattice reduct of a BL-algebra is distributive.
\end{Corollary}

BL-chains,  in turn,  admit a rich structure theory, as we proceed to explain. A \textit{t-norm} is a binary function $\ast \colon [0,1]^2\to[0,1]$ on the unit interval $[0, 1]$ that is commutative, associative, order preserving in both arguments, and such that $1 \ast a = a$, for every $a \in [0, 1]$. In addition, a t-norm is said to be \textit{continuous} when it is continuous with respect to the standard topology on $[0, 1]$. 
BL-chains are related to continuous t-norms as follows. On the one hand, every continuous t-norm $\ast$ induces a BL-chain
\[
\langle [0, 1]; \land, \lor, \ast, \to, 0, 1 \rangle,
\]
where $\land$ and $\lor$ are the binary operations of infinum and supremum with respect to the standard ordering of $[0, 1]$ and $\to$ is the binary operation defined by the rule
\[
a \to c \coloneqq \bigvee \{ b \in [0, 1] : b \ast a \leq c \}.
\]
BL-algebras of this form are known as \textit{standard}. On the other hand, every BL-chain embeds into an ultraproduct of standard BL-algebras \cite[Thm.\ 9]{Cignoli-Torrens:StandardCompletenessBL}.

In view of the theorem of Mostert and Shields \cite[Thm.\ B]{Mostert-Shields:OrdinalSum}, every continuous t-norm $\ast$ can be decomposed into an ordinal sum of three special t-norms: the truncated sum $a \ast_{\text{\emph \L}} c \coloneqq \max\{0,a+c-1\}$, the product $a \ast_{P}c \coloneqq ac$ and the minimum operation $a \ast_{G} c \coloneqq \min \{ a, c \}$.

Because of this, the standard BL-algebras $\standardL$, $\standardP$ and $\standardG$ induced, respectively, by the three basic continuous t-norms $\ast_{\text{\emph \L}}$, $\ast_{P}$, and $\ast_{G}$ stand out among BL-chains. Indeed, each of them induces a distinguished axiomatic extension of the basic logic $\textbf{BL}$. For instance, \textit{{\L}ukasiewicz logic} $\lulog$ is defined, for every set of formulas $\Gamma \cup \{ \varphi \}$, as
\[
\Gamma \vdash_{\lulog} \varphi \Longleftrightarrow \text{ there exists a finite }\Delta \subseteq \Gamma \text{ such that }\btau[\Delta] \vDash_{\standardL} \btau(\varphi),
\]
where $\btau \coloneqq \{ x \thickapprox 1 \}$. \textit{Product logic} $\prlog$ and \textit{G\"odel-Dummett logic} $\golog$ (sometimes called simply \textit{G\"odel logic}) are obtained similarly, replacing $\standardL$ by $\standardP$ and $\standardG$ respectively, see, e.g., \cite{CiMuOt99,Ha98}.

{\L}ukasiewicz, product, and G\"odel logic are algebraized, respectively, by varieties $\class{MV} \coloneqq \VVV(\standardL)$ of \textit{MV-algebras}, $\class{PA} \coloneqq \VVV(\standardP)$ of \textit{product algebras}, and $\class{GA} \coloneqq \VVV(\standardG)$ of \textit{G\"odel algebras}. Notably, 
\[
\class{MV} \coloneqq \QQQ(\standardL) \qquad\class{PA} \coloneqq \QQQ(\standardP) \qquad \class{GA} \coloneqq \QQQ(\standardG).
\]
The first equality above can be traced back to \cite[Lem.\ B]{Hay:Axiomatization}, see also the discussion in 
\cite{GpT98} or \cite[Lem.\ 3.2.11(3)]{Ha98}, 
the second is implicit in \cite{Hajek-Godo-Esteva:ProductLogic,Ha98} and is based on the fact that all nontrivial totally ordered abelian groups have the same universal theory  \cite{Gurevich-Kokorin:UnivEquivalence}, while the third is relatively straightforward.

Sufficiently well-structured MV-algebras, product algebras, and G\"odel algebras can be expanded with rational constants, as we proceed to explain. Consider a set of constants
\[
\mathcal{C}= \{\boldsymbol{c}_q : q \in [0, 1] \cap \mathbb{Q}\},
\] 
where $\mathbb{Q}$ denotes the set of rational numbers. Observe that $[0, 1] \cap \mathbb{Q}$ is the universe of a subalgebra of $\standardL$ (resp.\ of $\standardP$ and $\standardG$) that we denote by $\rationalL$ (resp.\  $\rationalP$ and $\rationalG$). Because of this, given an algebra $\A \in \{ \standardL, \standardP, \standardG \}$, we can consider the set $\mathcal{B}(\A)$ of equations in the language of $\A$ expanded with the constants in $\mathcal{C}$ of the form
\begin{align*}
\boldsymbol{c}_p \cdot \boldsymbol{c}_q \thickapprox \boldsymbol{c}_{p\cdot^{\A} q} \qquad \boldsymbol{c}_p \to \boldsymbol{c}_q \thickapprox \boldsymbol{c}_{p\to{^\A}q} \qquad \boldsymbol{c}_0 \thickapprox 0 \qquad \boldsymbol{c}_1 \thickapprox 1,
\end{align*}
for every $p, q \in \mathbb{Q} \cap [0, 1]$. The equations in $\mathcal{B}(\A)$  are sometimes called the \textit{bookkeeping axioms} of $\A$.  We do not include bookkeeping axioms for the lattice operations, because these can be defined in terms of $\cdot$ and $\to$.

\begin{law} An algebra $\A$ in the language of BL-algebras expanded with constants in $\mathcal{C}$ is said to be
\benroman
\item a \textit{rational MV-algebra} if the BL-reduct of $\A$ is an MV-algebra and $\A$ validates the bookkeeping axioms $\mathcal{B}(\standardL)$;
\item a \textit{rational product algebra} if the BL-reduct of $\A$ is a product algebra and $\A$ validates the bookkeeping axioms $\mathcal{B}(\standardP)$;
\item a \textit{rational G\"odel algebra} if the BL-reduct of $\A$ is a G\"odel algebra and $\A$ validates the bookkeeping axioms $\mathcal{B}(\standardG)$.
\eroman
We denote by $\class{RMV}, \class{RPA}$ and $\class{RGA}$ the varieties\footnote{Notice that $\class{RMV}, \class{RPA}$, and $\class{RGA}$ are varieties, because so are $\class{MV}, \class{PA}$, and $\class{GA}$, and the bookkeeping axioms are equations.} of rational MV-algebras, rational product algebras, and rational G\"odel algebras respectively.
\end{law}

Canonical rational MV, product, and G\"odel algebras can be obtained by expanding the standard BL-algebras $\standardL, \standardP$, and $\standardG$ with the natural interpretation of the constants in $\mathcal{C}$, that is, by interpreting $\boldsymbol{c}_q$ as the rational $q$. We denote these expansions, respectively, by $\standardLQ, \standardPQ$, and $\standardGQ$. Furthermore, we denote their subalgebras with universe $\mathbb{Q} \cap [0, 1]$ by  $\rationalLQ, \rationalPQ$, and $\rationalGQ$ respectively. The importance of the algebras $\standardLQ, \standardPQ$, and $\standardGQ$ is witnessed by the equalities 
\[
\class{RMV} = \VVV(\standardLQ) \qquad\class{RPA} = \VVV(\standardPQ)  \qquad\class{RGA} = \VVV(\standardGQ).
\]
For the second and the third equalities above, see \cite[Thm.\ 5.4]{Savicky-CEGN:ProductTruthConstants} and \cite[Thm.\ 13]{Esteva-Godo-Noguera:RationalWNM}.\ Notably, $\class{RMV}$ coincides also with the quasivariety generated by $\standardLQ$ \cite[Thm.\ 3.3.14]{Ha98}.  This contrasts with the case of $\class{RPA}$ and $\class{RGA}$, as there are not the quasivarieties generated by $\standardPQ$ and $\standardGQ$, see \cite[Lem.\ 3.6]{Savicky-CEGN:ProductTruthConstants} and  \cite[Sec.\ 4]{Esteva-Godo-Noguera:RationalWNM}.

From viewpoint of logic,  the varieties $\class{RMV}, \class{RPA}$, and $\class{RGA}$ algebraize expansions of $\lulog, \prlog$, and $\golog$. For instance $\class{RMV}$ algebraizes \textit{rational {\L}ukasiewicz logic} $\rlulog$ defined, for every set of formulas $\Gamma \cup \{\varphi \}$ as
\[
\Gamma \vdash_{\rlulog} \varphi \Longleftrightarrow \text{ there exists a finite }\Delta \subseteq \Gamma \text{ such that }\btau[\Delta] \vDash_{\class{RMV}} \btau(\varphi),
\]
where $\btau \coloneqq \{ x \thickapprox 1 \}$.  \textit{Rational product logic} $\rprlog$ and \textit{rational G\"odel logic} $\rgolog$ are obtained similarly, replacing $\class{RMV}$ by $\class{RPA}$ and $\class{RGA}$, see, e.g., \cite{Esteva-Godo-Noguera:Rational}.

Notice that Theorem \ref{Thm: filters and congruences} and Lemma \ref{lem:: prime filter lemma} apply to rational {\L}ukasiewicz,  rational product, and rational G\"odel algebras as well, because the addition of constants to a given algebra does not change its congruences and filters.

\section{Extensions of rational product logic}\label{Sec:5}

In view of the following result, the lattice of extensions of product logic (without rational constants) is a three-element chain:

\begin{Theorem}[\protect{\cite[Cor.\ 3.22]{CiMe09}}]\label{Thm:no-constants-product}
The unique proper nontrivial extension of $\prlog$ is algebraized by a variety term-equivalent to that of Boolean algebras. Consequently, every extension of $\prlog$ is axiomatic and $\prlog$ is HSC.
\end{Theorem}

It is therefore sensible to wonder whether a similar situation holds for the case of rational product logic $\rprlog$. The aim of this section is to shed light on this problem. To this end, it is convenient to separate the case of axiomatic extensions from that of arbitrary extensions of $\rprlog$. This is because, as regarding axiomatic extensions, $\rprlog$ behaves similarly to $\prlog$.

\begin{Theorem}\label{Thm:rational-variety}
There are only two proper nontrivial subvarieties $\class{K}_{1}$ and $\class{K}_{2}$ of $\class{RPA}$. 
\benroman
\item $\class{K}_{1}$ is term-equivalent to the variety of Boolean algebras and is axiomatized by the pair of equations $\boldsymbol{c}_{q} \thickapprox 1$ and $x \lor (x \to 0) \thickapprox 1$, where $q$ is any rational number in the interval $(0,1)$;
\item $\class{K}_{2}$ is term-equivalent to the variety of product algebras and is axiomatized by the equation $\boldsymbol{c}_{q} \thickapprox 1$, where $q$ is any rational number in the interval $(0,1)$.
\eroman
Consequently, $\class{K}_{1} \subsetneq \class{K}_{2}$ and $\mathcal{V}(\class{RPA})$ is a four-element chain.
\end{Theorem}

Since the lattice of axiomatic extensions of $\rprlog$ is dually isomorphic to that of subvarieties of $\class{RPA}$, the above result can be rephrased in logical parlance as follows.

\begin{Corollary}
The lattice of axiomatic extensions of $\rprlog$ is a four-element chain. The two sole proper consistent axiomatic extensions of $\rprlog$ are algebraized by varieties term-equivalent, respectively, to those of Boolean and product algebras.
\end{Corollary}

On the other hand, the lattice of extensions of rational product logic is quite complicated.

\begin{Theorem}\label{Thm:Q-universal}
The variety of rational product algebras is $\mathcal{Q}$-universal. Consequently, the lattice of extensions of $\rprlog$ has the cardinality of the continuum and does not validate any nontrivial lattice equation.
\end{Theorem}

\noindent Accordingly, from the point of view of extensions, $\rprlog$ is by far richer than $\prlog$.

The remaining part of this section is devoted to the proofs of Theorems \ref{Thm:rational-variety} and \ref{Thm:Q-universal}. In order to establish Theorem \ref{Thm:rational-variety}, we rely on the following observation.

\begin{Proposition}\label{Prop:Zuzka-trick}
$\class{RPA} = \VVV(\rationalPQ)$. 
\end{Proposition}

\begin{proof}
In \cite[Thm.\ 5.4]{Savicky-CEGN:ProductTruthConstants} it is shown that $\class{RPA} = \VVV(\standardPQ)$. Accordingly, to prove that $\class{RPA} = \VVV(\rationalPQ)$, it suffices to show that if an equation fails in $\standardPQ$, then it also fails in $\rationalPQ$. 
Assume, towards a contradiction, that there is an equation $\epsilon(\vec x)\thickapprox\delta(\vec x)$ 
true in $\rationalPQ$ and a tuple $\vec a$ of reals in $[0,1]$ 
such that $\epsilon^{\standardPQ}(\vec a)\neq\delta^\standardPQ (\vec a)$. 
First notice that we can assume $\delta$ to be $1$; otherwise we replace the original equation with 
$\epsilon \leftrightarrow \delta \thickapprox 1$. 
In the rest of the proof we moreover assume that the lattice connectives do not occur in $\epsilon$;
this is without loss of generality as they are term-definable from $\cdot$ and $\to$.

Recall that for $a,b\in [0,1]$,  we have
\begin{align*}
a \cdot^{\standardPQ} b = 0 &\Longleftrightarrow 0 = \min \{ a, b \}\\
a \to^{\standardPQ} b = 0 &\Longleftrightarrow a > b = 0.
\end{align*}
Clearly,  $F \coloneqq(0, 1]$ is a filter on $\standardPQ$ and $\standardPQ/F$ is a two-element chain;
we can assume that its universe is $\{ 0, 1 \}$ (under the identification of $0$ with $0/ F$ and $1$ with 
$1 / F$). For $a,b\in [0,1]$, we write $a\sim b$ if and only if $a/F = b/F$. 

% Let $\phi$ be the congruence of $\standardPQ$ with blocks $\{ 0 \}$ and $(0, 1]$, 
% and set $\A \coloneqq \standardPQ / \phi$. As $\A$ is a two-element algebra, 
% we can assume that its universe is $\{ 0, 1 \}$ (under the identification of 
% $0$ with $0/ \phi$ and $1$ with $1 / \phi$). 
% Let also $h \colon \standardPQ \to \A$ be the canonical surjection, 
% which sends all nonzero elements to $1$. For any RPA-term $\varphi(\vec{x})$ 
% and tuple $\vec{a} \in [0, 1]$, one can label each subterm $\psi$ with $h(\psi^{\standardPQ}(\vec{a}))$, 
% the latter indicating whether $\psi^{\standardPQ}(\vec{a})$ is $0$. 
% Notice that this labelling is invariant under any perturbation of $\vec{a}$ 
% that preserves $h$ on the values of the variables $\vec{x}$.

%Let $\standardPQpower{+}$ be the zero-free subreduct of $\standardPQ$ with the universe $(0, 1]$. 
%In the algebra $\standardPQpower{+}$, both $\cdot$ and $\to$ are interpreted with continuous functions, so any   zero-free term defines a continuous function in $\standardPQpower{+}$ and its interpretation is hence fully determined by the values assigned to the rationals. In particular, for any $a \in (0,1]$, a sequence 
%$\{a_k\}_{k\in \omega}$ of rationals $0<a_k\leq 1$ tending to $a$, and a zero-free term $\varphi(x,\vec{y})$, the values $\{ \varphi^{\standardPQpower{+}}(a_k,\vec{c}) : k \in \omega \}$ tend to $\varphi^{\standardPQpower{+}}(a,\vec{c})$ for any choice of $\vec{c}$ in $(0, 1]$.   

Now recall that  $\epsilon^{\standardPQ}(\vec a)<1$ for some $\vec a =\langle a_{1}, \dots, a_{n} \rangle$ in $[0,1]^n$.
 Then take $I=\{i \leq n :  a_i \mbox{ is irrational}\}$. For an arbitrarily chosen $i\in I$, 
fix a sequence $\{ a_{ik} : k \in \omega\}$ of rationals in $(0, 1]$ tending to $a_i$ 
(this is possible because $a_{i}$, being irrational, is positive). 
To conclude the proof, it is enough to find, for the chosen $i \in I$, a rational $c_{i} \in [0, 1]$ such that
\[
\epsilon^{\standardPQ}(a_{1}, \dots, a_{i-1}, c_{i}, a_{i+1}, \dots, a_{n}) < 1.
\]
This is sufficient as the process can be iterated for the remaining elements of $I\smallsetminus\{ i\}$, 
finally obtaining  rationals $c_{1}, \dots, c_{n} \in [0, 1]$ such that $\epsilon^{\standardPQ}(c_{1}, \dots, c_{n}) < 1$, as desired.

Accordingly, fix an $i \in I$ and for each subterm $\eta(\vec x)$ of $\epsilon(\vec x)$,
 let $f_{\eta} \colon (0, 1] \to [0, 1]$ be the map defined by the rule
\[
f_{\eta}(z) \coloneqq \eta^{\standardPQ}(a_{1}, \dots, a_{i-1}, z, a_{i+1}, \dots, a_{n}).
\]
We claim that the sequence $\{ f_{\eta}(a_{ik}) : k \in \omega \}$ tends to $\eta^{\standardPQ}(a_{1}, \dots, a_{n})$ for any choice of $\eta$ a subterm of $\epsilon$. The proof is by induction on term structure of $\eta$. 
The cases where $\eta$ is a variable or a constant are straightforward. 
Given that $a_{i}$ and all $a_{ik}$ are positive, we have $a_i \sim a_{ik}$ for each $k\in\omega$.
Then for every subterm $\eta$ of $\epsilon$ and every $k \in \omega$,
\[
f_{\eta}(a_{i}) \sim f_{\eta}(a_{ik}). 
\]
Consequently, if $\eta^{\standardPQ}(a_{1}, \dots, a_{n}) = 0$, then $\{ f_{\eta}(a_{ik}) : k \in \omega \}$ is a constant sequence of zeros and we are done. Then we consider the case where $\eta^{\standardPQ}(a_{1}, \dots, a_{n}) \ne 0$. 

For the inductive step, observe that if
$\eta$ is of the form $\varphi_{1} \cdot \varphi_{2}$, the result follows from the inductive hypothesis and the fact that $\cdot$ is continuous in $[0, 1]$. Then we consider the case where $\eta$ is of the form $\varphi_{1} \to \varphi_{2}$. We have
\[
0 \ne \eta^{\standardPQ}(\vec{a}) = \varphi_{1}^{\standardPQ}(\vec{a}) \to^{\standardPQ} \varphi_{2}^{\standardPQ}(\vec{a}).
\]
Consequently, either $\varphi_{1}^{\standardPQ}(\vec{a}) = 0$ or $\varphi_{1}^{\standardPQ}(\vec{a}), \varphi_{2}^{\standardPQ}(\vec{a}) > 0$. Suppose that the latter holds. Then $f_{\varphi_j}(a_{ik}) \sim f_{\varphi_j}(a_i) = \varphi_{j}^{\standardPQ}(\vec{a})>0$, 
for all $j = 1, 2$ and $k \in \omega$. Hence, 
$0 < f_{\varphi_j} (a_{ik})$ and the result follows from the inductive hypothesis and the fact that $\to$ is continuous in $(0, 1]$. It only remains to consider the case where $\varphi_{1}^{\standardPQ}(\vec{a}) = 0$.   We have $\varphi_{1}^{\standardPQ}(\vec{a}) \to \varphi_{2}^{\standardPQ}(\vec{a}) = 1$. 
Furthermore, as $a_i \sim a_{ik}$ for $k\in \omega$, we have 
 $\varphi_{1}^{\standardPQ}(a_{1}, \dots, a_{i-1}, a_{ik}, a_{i+1}, \dots, a_{n}) = 0$ for all $k \in \omega$, whence
\begin{align*}
& \,\varphi_{1}^{\standardPQ}(a_{1}, \dots, a_{i-1}, a_{ik}, a_{i+1}, \dots, a_{n}) \to^{\standardPQ} \varphi_{2}^{\standardPQ}(a_{1}, \dots, a_{i-1}, a_{ik}, a_{i+1}, \dots, a_{n})\\
=& \,0 \to^{\standardPQ} \varphi_{2}^{\standardPQ}(a_{1}, \dots, a_{i-1}, a_{ik}, a_{i+1}, \dots, a_{n})\\
=& \, 1 = \varphi_{1}^{\standardPQ}(\vec{a}) \to^{\standardPQ} \varphi_{2}^{\standardPQ}(\vec{a}) \text{, for all }k \in \omega.
\end{align*}
Thus, $1=f_\eta(a_i)=f_{\eta}(a_{ik})$ for each $k\in \omega$. This establishes the claim.

Given the claim, the fact that $\epsilon^{\standardPQ}(a_{1}, \dots, a_{n}) < 1$ implies the existence of an $m \in \omega$ such that $\epsilon^{\standardPQ}(a_{1}, \dots, a_{i-1}, a_{im}, a_{i+1}, \dots, a_{n})<1$. 
Taking $c_{i} \coloneqq a_{im}$, we are done.
\end{proof}

\begin{proof}[Proof of Theorem \ref{Thm:rational-variety}.]
Observe that a rational product algebra validates the equation $\boldsymbol{c}_{q} \thickapprox 1$ for some $q \in (0, 1) \cap \mathbb{Q}$ if and only if it validates all the equations  $\{ \boldsymbol{c}_{p}\thickapprox 1 : p \in (0, 1) \cap \mathbb{Q}\}$. 
Indeed, let $p,q \in (0, 1) \cap \mathbb{Q}$ and $c_q^\A=1^\A$. Let $n$ be an integer such that $q^n\leq p$. Then, by the bookkeeping axioms, we have
\[
c_p^\A\geq c_{q^n}^\A= (c_q^\A)^n=1^\A.
\]

Let $\class{PA}^{\ast}$ be the class of rational product algebras $\A$ in which $c_{q}^{\A} = 1$ for all rational numbers $q \in (0, 1]$. Clearly, $\class{PA}^{\ast}$ is a subvariety of $\class{RPA}$ term-equivalent to that of product algebras. Thus, in view of Theorem \ref{Thm:no-constants-product}, it suffices to prove that $\class{PA}^{\ast}$ is the largest proper subvariety of $\class{RPA}$. To this end, let $\A$ be a rational product algebra such that $\VVV(\A)$ is a proper subvariety of $\class{RPA}$. By Proposition \ref{Prop:Zuzka-trick}, the zero-generated subalgebra $\C$ of $\A$ cannot be isomorphic to $\rationalPQ$. Accordingly, there is a filter $F$ of $\rationalPQ$, different from $\{1\}$, such that the algebra $\C$ is isomorphic to $\rationalPQ / F$. The algebra $\rationalPQ$ has only three filters: $\{1\}$, $(0,1]$ and $[0,1]$. Thus $(0,1]\subseteq F$, and hence $c_{q}^{\A} = 1$ for all $q \in (0, 1]$. This means that $\A \in \class{PA}^{\ast}$. This shows that  $\class{PA}^{\ast}$ is the largest proper subvariety of $\class{RPA}$.
\end{proof}

The remaining part of the section is devoted to the proof of Theorem \ref{Thm:Q-universal}. 
%However, this condition is not quite suited for our need.\ 
%Indeed, it is assumed there that the language is finite, which clearly fails for $\class{RPA}$, and that the algebras witnessing this condition are finite, which cannot be achieved here since  every finite rational product algebra is term-equivalent to a Boolean algebra.
Let $\Prime$ be the set of prime numbers and let $\mathcal{P}_{< \omega}(\Prime)$ be the set of finite subsets of $\Prime$. We denote by $\mathcal{S}(\mathcal{P}_{<\omega}(\Prime))$ the lattice of universes of subalgebras of $\langle \mathcal{P}_{< \omega}(\Prime), \cup, \emptyset \rangle$ with set inclusion as the order. The following observation can be extracted from the proof of \cite[Thm.\ 3.3]{AdaDzi94}.

\begin{Theorem}
\label{Theorem:: Adams-Dziobiak}
Let  ${\class K}$ be a quasivariety. If there exist a subquasivariety $\class{M}$ of $\class{K}$ and a surjective bounded-lattice homomorphism
\[
h \colon \mathcal{Q}(\mathsf{M}) \to \mathcal{S}(\mathcal{P}_{<\omega}(\Prime)),
\]
then $\class{K}$ is $\mathcal Q$-universal.
\end{Theorem}

\begin{proof}[Proof of Theorem \ref{Thm:Q-universal}]
For the variety $\class{RPA}$ we find a quasivariety $\class{M}$ and a surjective homomorphism $h$ as in Theorem \ref{Theorem:: Adams-Dziobiak}.

We begin by defining a family of rational product algebras $\{ \A_X : X \in \mathcal{P}_{<\omega}(\Prime) \}$. For a finite nonempty set $X$ of prime numbers, let $\A_X$ be the subalgebra of $\standardPQpower{X}$ generated by the function  $\mathrm{inv}_X\colon X\to [0,1]$ that sends an element $p \in X$ to $\nicefrac{1}{\sqrt{p}}$. Notice that $\A_\emptyset$ is a trivial algebra.
Let
\[
\class{M}\coloneqq \QQQ(\{\A_X : X\in \mathcal{P}_{<\omega}(\Prime)\})
\]
and, for a subquasivariety $\class{N}$ of $\mathsf{M}$, define
\[
h(\class{N})\coloneqq \{X\in \mathcal{P}_{<\omega}({\Prime}) : \A_{X} \in\class{N} \}.
\]

In order to prove that $h$ is a well-defined map from $\mathcal{Q}(\mathsf{M})$ to $\mathcal{S}(\mathcal{P}_{<\omega}(\Prime))$, let $\class{N}\in\mathcal{Q}(\mathsf{M})$.
Since $\A_\emptyset$ is a trivial algebra, $\A_\emptyset\in\class{N}$, and hence $\emptyset\in h(\class{N})$. To prove that $h(\class{N})$ is closed under binary unions, consider $X_1, X_2 \in h(\class{N})$. By definition of $h$, we have $\A_{X_1}, \A_{X_2} \in \class{N}$. Moreover, the map that sends an element $a \in A_{X_1\cup X_2}$ to the pair $(a{\upharpoonright}_{X_1},a{\upharpoonright}_{X_2})$ is an embedding of $\A_{X_1\cup X_2}$ into $\A_{X_1}\times\A_{X_2}$. Hence,
\[
\A_{X_1\cup X_2} \in \III\SSS(\A_{X_1} \times \A_{X_2}) \subseteq \III\SSS\PPP(\class{N}) \subseteq \class{N}.
\]
As a consequence, $X_1 \cup X_2 \in h(\class{N})$. We conclude that $h \colon \mathcal{Q}(\mathsf{M}) \to \mathcal{S}(\mathcal{P}_{<\omega}(\Prime))$ is well defined, as desired.

It follows from the definition of $h$ that it preserves the binary meet. We also have $h(\class{M})=\mathcal{P}_{<\omega}(\Prime)$, which means that $h$ sends the top element of $\mathcal{Q}(\mathsf{M})$ to the top element of $\mathcal{S}(\mathcal{P}_{<\omega}(\Prime))$. Let $\class{Tr}$ be the class of all trivial algebras (in the language of $\class{RPA}$). Since $\A_\emptyset$ is the only trivial algebra in the family $\{\A_X : X \in\mathcal{P}_{<\omega}(\Prime)\}$, we have $h(\class{Tr})=\{\emptyset\}$. This means that $h$ sends the bottom element of $\mathcal{Q}(\mathsf{M})$ to the bottom element of $\mathcal{S}(\mathcal{P}_{<\omega}(\Prime))$. 

Thus, it only remains to show that $f$ preserves the binary join and is surjective. The proof will proceed through a series of claims. For every $X\in \mathcal{P}_{<\omega}(\Prime)$, 
let us consider the formulas
\[
\Gamma_X(z)\coloneqq  \bigvee_{p\in X} (z^2\leftrightarrow \boldsymbol{c}_{\nicefrac{1}{p}}) 
\thickapprox 1 
\]
 and 
\[
\Delta_X(z)\coloneqq \Gamma_X(z) \boland \foo_{p\in X}\mathrel{\boneg} \Gamma_{X\smallsetminus \{p\}}(z). 
\]
We interpret $\Gamma_\emptyset(z)$ and $ \Delta_\emptyset(z)$ as $0\thickapprox 1$. Hence these formulas hold only in trivial algebras. Moreover, for every $X\in\mathcal{P}_{<\omega}(\Prime)$ and $a \in [0, 1]$,
\[
\standardPQ\models \Gamma_X(a) \Longleftrightarrow a\in\{\nicefrac{1}{\sqrt{p}} : p\in X\}.
\]
As a consequence, we obtain $\A_X\models\Delta_X(\mathrm{inv}_X)$.

For every $p \in \Prime$, let $\A_p$ be a subalgebra of $\standardPQ$ generated by $\nicefrac{1}{\sqrt{p}}$ (notice that $\A_p\cong\A_{\{p\}}$).

\begin{Claim}
\label{Claim: 1/sqrt{q} not in Ap}
Let $p$ and $q$ be distinct prime numbers. Then $\nicefrac{1}{\sqrt{p}}$ does not belong to $\A_q$.
\end{Claim}

\begin{proof}[Proof of the claim.]
Let 
\[
B_q\coloneqq \left(\Rat\cup \{a\cdot\sqrt{q} : a\in \Rat\}\cup \{{\large\nicefrac{a}{\sqrt{q}}} : a\in \Rat\}\right)\cap [0,1].
\]
Then $B_q$ is the universe of a subalgebra $\B_q$ of $\standardPQ$. Since $\nicefrac{1}{\sqrt{q}}\in B_{q}$, we have $A_{q}\subseteq B_{q}$. And hence, since $\nicefrac{1}{\sqrt{p}}\not\in B_{q}$, we have $\nicefrac{1}{\sqrt{p}}\not\in A_{q}$. 
\end{proof}

\begin{Claim}
\label{Claim:: A_Y models Delta_X => X=Y}
Let $X,Y\in \mathcal{P}_{<\omega}(\Prime)$ and $a\in A_Y$. Then $\A_Y\models\Delta_X(a)$ if and only if  $X=Y$ and $a=\mathrm{inv}_X$.
\end{Claim}

\begin{proof}[Proof of the claim.]
The implication from right to left holds because $\A_X\models\Delta_X(\mathrm{inv}_X)$, for all $X \in \mathcal{P}_{<\omega}({\Prime})$. To prove the other implication, suppose that $\A_Y\models\Delta_X(a)$. Recall that $\A_Y$ is the subalgebra of $\standardPQpower{Y}$ generated by $\mathrm{inv}_Y$. Therefore, $a$ is a function $a \colon Y \to [0, 1]$. Let $q\in Y$. Since $\A_q$ is a homomorphic image of $\A_Y$ and  $\Gamma_X(z)$ is an equation such that $\A_Y\models\Delta_X(a)$, we have $\A_q\models\Gamma_X(a(q))$. As $\A_q$ is a subalgebra of $\standardPQ$, we obtain that $a(q)\in\{\nicefrac{1}{\sqrt{p}} : p\in X\}\cap A_q$.  By Claim \ref{Claim: 1/sqrt{q} not in Ap}, it follows that $a(q)=\nicefrac{1}{\sqrt{q}}$ and $q\in X$. This shows that $a=\mathrm{inv}_Y$ and $Y\subseteq X$. It only remains to show that $X \subseteq Y$. Accordingly, let $p\in X$.  Since $\A_Y \vDash \Delta_X(a)$ and $a=\mathrm{inv}_Y$, we have $\A_Y\models\Gamma_Z(\mathrm{inv}_Y)$ for every finite set $Z$ of prime numbers such that $Y\subseteq Z$ and $\A_Y\nvDash\Gamma_{X\smallsetminus\{p\}}(\mathrm{inv}_Y)$. Consequently, $Y\nsubseteq X\smallsetminus\{p\}$. As $Y \subseteq X$, this implies $p\in Y$. Hence, we conclude that $X \subseteq Y$, as desired.
\end{proof}

\begin{Claim}
\label{Claim:: dealing with products}
Let $X\in \mathcal{P}_{<\omega}(\Prime)$ and let $\{ \B_i : i \in I \}$ be a family of rational product algebras, each with an element $b_i \in B_i$. Furthermore, let $\B = \prod_{i\in I}\B_i$ and let $b$ be the element of $B$ whose $i$-th coordinate is $b_i$. If  $\B\models \Delta_X(b)$, then there exists a family of sets $\{ Y_i : i \in I \}$ such that
\benroman
\item $X=\bigcup_{i\in I}Y_i$ and
\item $\B_i\models\Delta_{Y_i}(b_i)$, for every $i\in I$.
\eroman
\end{Claim}

\begin{proof}[Proof of the claim.]
Let $i\in I$. Since $\B\models \Delta_X(b)$, we have $\B\models \Gamma_X(b)$. As $\Gamma_X(z)$ is an equation, it is preserved by homomorphisms and, therefore, $\B_i\models\Gamma_X(b_i)$. Hence, since $X$ is finite, there exists  a subset $Y_i$ of $X$  such that $\B_i\models\Gamma_{Y_i}(b_i)$ but $\B_i\not\models\Gamma_{Y_i\smallsetminus\{p\}}(b_i)$ for every $p\in Y_i$. In particular, $\B_i\models\Delta_{Y_i}(b_i)$.

Let then $Y\coloneqq\bigcup_{i\in I}Y_i$. By construction, $Y\subseteq X$. In order to prove reverse inclusion, consider $i \in I$. From $Y_i \subseteq Y$ and $\B_i \models \Gamma_{Y_i}(b_i)$ it follows $\B_i\models \Gamma_Y(b_i)$. Thus, $\B\models\Gamma_Y(b)$ and, therefore, $\B\models\Gamma_Z(b)$ for every finite set $Z \supseteq Y$ of primes. Let $p\in X$. Since $\B \vDash \Delta_X(b)$, we have $\B\nvDash\Gamma_{X\smallsetminus\{p\}}(b)$. It follows that $Y\not\subseteq X\smallsetminus\{p\}$.  As $Y \subseteq X$, this implies $p \in Y$. Hence, we conclude that  $X\subseteq Y$. 
\end{proof}

We are now in a position to prove surjectivity.

\begin{Claim}
The map $h$ is surjective.
\end{Claim}

\begin{proof}[Proof of the claim.]
Let $S\in \mathcal{S}(\mathcal{P}_{<\omega}(\Prime))$ and define
\[
\class{N}\coloneqq \QQQ(\{\A_Y : Y\in S\}).
\] 
We clearly have $S\subseteq h(\class{N})$. For the verification of the reverse inclusion, let us consider a set $X$ in $h(\class{N})$, i.e., such that $\A_X\in \class{N}$.  Then
\[
\A_X \in \QQQ(\{\A_Y : Y\in S\}) = \III\SSS\PPP\PPU(\{\A_Y : Y\in S\}).
\]
Accordingly, there is a family of algebras $\{\B_i : i \in I \}$ such that $\A_X$ embeds into the product $\prod_{i\in I} \B_i$ and, for every $i \in I$, there is a family of algebras $\{ \C_{ij} : j \in J_i \}$ and an ultrafilter $U_i$ on $J_i$ such that $\B_i=\prod_{j\in J_i}\C_{ij}/U_i$ and the various $\C_{ij}$ belong to the family $\{\A_Y : Y\in S\}$.
 
Since $\A_X\models\Delta_X(\mathrm{inv}_X)$  and $\A_X$ embeds into $\prod_{i \in I}\B_i$,  there is $b \in B$ such that $\prod_{i\in I}\B_i\models\Delta_X(b)$. Then we may apply Claim  \ref{Claim:: dealing with products} obtaining a family $\{ Y_i : i \in I \}$ of finite subsets of $X$ as in the statement of the claim. In particular, for every $i \in I$, we have  $\B_i\models\exists z\,\Delta_{Y_i}(z)$. By \L{}o\'s' Theorem \cite[Thm. V.2.9]{BuSa00}, there exists a nonempty set $W\in U_i$ such that $\C_{ij}\models \exists z\,\Delta_{Y_i}(z)$ for every $j\in W$. Consequently, there exists $j\in J_i$ such that
$\C_{ij}\models \exists z\,\Delta_{Y_i}(z)$. Recall that $\C_{ij}=\A_Y$ for some $Y\in S$. By Claim \ref{Claim:: A_Y models Delta_X => X=Y}, the sets $Y_i$ and $Y$ are equal. In this way, we showed that $Y_i\in S$, for every $i\in I$. Since $X$ is finite, the set $\{Y_i : i\in I\}$, consisting of subsets of $X$, is finite as well. Moreover, it is a subset of $S$. Thus, its union, which equals $X$, belongs to  $S$. 
\end{proof}

\begin{Claim}
\label{Claim:: dealing with ultraproducts}
Let $X\in \mathcal{P}_{<\omega}(\Prime)$. If $\B\models\exists z\, \Delta_X(z)$ and $\B \in \PPU(\{\A_Y : Y\in \mathcal{P}_{<\omega}(\Prime)\})$, then $\B \in \III\PPU(\A_X)$.
\end{Claim}

\begin{proof}[Proof of the claim.]
Assume that $\B\models\exists z\, \Delta_X(z)$ and that $\B=\prod_{i\in I} \C_i/U$, where all $\C_i$ belong to $\{\A_Y : Y\in \mathcal{P}_{<\omega}(\Prime)\}$.\ 
Let $J\coloneqq\{ i \in I : \C_i\models \exists z\, \Delta_X(z)\}$. By 
\L{}o\'s' Theorem, $J\in U$. Moreover, by Claim \ref{Claim:: A_Y models Delta_X => X=Y}, $\C_i=\A_X$ for every $i\in J$. Thus $\B$ is isomorphic to the ultrapower $\A_X^J/(U\cap\mathcal P(J))$ of $\A_X$.
\end{proof}

\begin{Claim}
\label{Claim:: A_X is a subalgebra}
Let $X\in \mathcal{P}_{<\omega}(\Prime)$ and $\A\in \class{M}$. If $\A\models\exists z\,\Delta_X(z)$, then $\A_X$ embeds into $\A$. 
\end{Claim}

\begin{proof}[Proof of the claim.]
Let $a\in \A$ be such that $\A\models\Delta_X(a)$. We will show that the subalgebra $\A'$ of $\A$ generated by $a$ is isomorphic to $\A_X$.

As $\A \in \class{M} =\III\SSS\PPP\PPU(\{\A_Y : Y\in S\})$, there exists an embedding $e\colon \A\to\prod_{i\in I} \B_i$, where 
 $\B_i=\prod_{j\in J_i}\C_{ij}/U_i$ for some algebras  $\C_{ij}\in \{\A_Y : Y\in \mathcal{P}_{<\omega}(\Prime)\}$ and ultrafilters $U_i$ on $J_i$.  Let $b \coloneqq e(a)$. Then, since $e$ is an embedding, $\prod_{i\in I} \B_i \models \Delta_X(b)$ and, hence, we may apply Claim \ref{Claim:: dealing with products} obtaining a family $\{ Y_i : i \in I \}$ of finite subsets of $X$ as in the statement of the claim. Consider $i \in I$. By Claim \ref{Claim:: dealing with products},  $\B_i\models \Delta_{Y_i}(b_i)$,    where $b_i$ is the $i$-th coordinate of $b$. Moreover, by Claim \ref{Claim:: dealing with ultraproducts}, $\B_i$ is isomorphic to an ultrapower of $\A_{Y_i}$. In particular, there exists an elementary embedding $e_i\colon \A_{Y_i}\to\B_i$.   Since $\A_{Y_i}\models \Delta_{Y_i}(\mathrm{inv}_{Y_i})$ and $e_i$ is an embedding, $\B_i\models \Delta_{Y_i}(e_i(\mathrm{inv}_{Y_i}))$. Furthermore, by Claim
\ref{Claim:: A_Y models Delta_X => X=Y}, there is at most one element of $\A_{Y_i}$ satisfying $\Delta_{Y_i}(z)$. Since $\A_{Y_i}$ are $\B_i$ are elementarily equivalent, the same is true for $\B_i$. Therefore, from  $\B_i\models \Delta_{Y_i}(e_i(\mathrm{inv}_{Y_i})) \boland \Delta_{Y_i}(b_i)$ it follows $b_i=e_i(\mathrm{inv}_{Y_i})$.  
 
 Let $\B_i'$ be the subalgebra of $\B_i$ generated by $b_i$.   Since $\A_{Y_i}$ is generated by  $\mathrm{inv}_{Y_i}$ and $e_i \colon \A_{Y_i} \to \B_i$ is an embedding that sends $\mathrm{inv}_{Y_i}$ to $b_i$, the algebras $\A_{Y_i}$ and $\B_i'$ are isomorphic under $e_i\colon \A_{Y_i}\to\B'_i$.  Furthermore, recall that $\A'$ is the subalgebra of $\A$ generated by $a$ and that the map $e \colon \A \to \prod_{i \in I}\B_i$ is an embedding that sends $a$ to $b$. Since $b_i \in B_i'$ for all $i \in I$, the map $e$ restricts to an embedding of $\A'$ into $\prod_{i \in I}\B_i'$. Lastly, as each $e_i^{-1} \colon\B_i^{'} \to \A_{Y_i}$ is an isomorphism that sends $b_i$ to $\mathrm{inv}_{Y_i}$, there exists an embedding $e^+ \colon \A' \to \prod_{i \in I}\A_{Y_i}$ such that
 \[
 e^+(a) = \langle \mathrm{inv}_{Y_i} : i \in I \rangle.
 \]
Since $\bigcup_{i\in I}Y_i=X$, the subalgebra of $\prod_{i\in I}\A_{Y_i}$ generated by $e^+(a)$ is isomorphic to $\A_X$. Thus, $\A'$ is isomorphic to $\A_X$, as desired.
\end{proof}

\begin{Claim}
The map $h$ preserves binary joins.
\end{Claim}

\begin{proof}[Proof of the claim.]
Let $\class{N}_1$ and $\class{N}_2$ be subquasivarieties of $\class{M}$. We have to show that $h(\class{N}_1\lor \class{N}_2)=h(\class{N}_1)\lor h(\class{N}_2)$. Since $h$ is order preserving, the inclusion $h(\class{N}_1)\lor h(\class{N}_2)\subseteq h(\class{N}_1\lor \class{N}_2)$ holds. 

In order to verify the reverse inclusion, we consider a set  $X$ in $h(\class{N}_1\lor \class{N}_2)$. Then $\A_X\in \class{N}_1\lor \class{N}_2=\III\SSS\PPP\PPU(\class{N}_1\cup \class{N}_2)$.\  It is a general fact about ultraproducts that $\PPU(\class{N}_1 \cup \class{N}_2) = \PPU(\class{N}_1) \cup \PPU(\class{N}_2)$. Since $\class{N}_1$ and $\class{N}_2$ are closed under the formation of ultraproducts, this yields $\PPU(\class{N}_1\cup \class{N}_2)=\class{N}_1\cup \class{N}_2$. Since they are also closed under direct proucts,  there are $\B_1\in\class{N}_1$ and $\B_2\in\class{N}_2$ such that $\A_X$ embeds into into $\B_1\times \B_2$. Since $\A_X\models\exists z\, \Delta_X(z)$, we obtain $\B_1 \times \B_2\models\exists z\, \Delta_X(z)$. By Claim \ref{Claim:: dealing with products}, there are sets $Y_1$ and $Y_2$ such that $X=Y_1\cup Y_2$ and $\B_1\models\exists z\, \Delta_{Y_1}(z)$ and $\B_2\models\exists z\, \Delta_{Y_2}(z)$. By Claim \ref{Claim:: A_X is a subalgebra}, $\A_{Y_1}$ embeds into $\B_1$. Hence, $\A_{Y_1}\in\class{N}_1$ and, therefore, $Y_1\in h(\class{N_1})$. In the same way, we obtain that $Y_2\in h(\class{N_2})$. Thus $X=Y_1\cup Y_2\in h(\class{N_1})\lor h(\class{N_2})$.
\end{proof}

Hence, we conclude that $h \colon \mathcal{Q}(\mathsf{M}) \to \mathcal{S}(\mathcal{P}_{<\omega}(\Prime))$ is a surjective bounded-lattice homomorphism.
\end{proof}

\section{Structural completeness in rational product logic}\label{Sec:6}

While it is well known that product logic $\prlog$ is HSC (Theorem \ref{Thm:no-constants-product}), from Theorem \ref{Thm:Q-universal} and Corollary \ref{Cor:HSC-distributive} it follows that $\rprlog$ is not HSC. Indeed, as we shall see, $\rprlog$ is not even passively structurally complete. The next result provides an answer to the question which rules are admissible in $\rprlog$.

\begin{Theorem}\label{Thm:base-product}
The structural completion of\/ $\rprlog$ is the unique extension of $\rprlog$ whose equivalent algebraic semantics is $\QQQ(\rationalPQ)$. A base for the admissible rules of\/ $\rprlog$ is given by the set of rules of the form
\[
\boldsymbol{c}_{q} \lor z \rhd z \quad \quad (\boldsymbol{c}_{p} \leftrightarrow x^{n}) \lor z \rhd z,
\]
for each (equiv.\ some) $q \in (0, 1) \cap \mathbb{Q}$ and each $p \in [0, 1]\cap \mathbb{Q}$, $n \in \omega$ such that $\sqrt[n]{p}$ is irrational.
\end{Theorem}

The core of the proof of Theorem \ref{Thm:base-product} amounts to the following description of the universal theory of $\rationalPQ$.

\begin{Theorem}\label{Thm:axiom-chain-Q}
The universal theory of $\rationalPQ$ is axiomatized relative to $\class{RPA}$ by the sentences
\[
\forall xy\, ( x \leq y \boor y \leq x),\quad \boldsymbol{c}_{q} \not\thickapprox 1, \, \, \text{ and } \, \,  \forall x\,(\boldsymbol{c}_{p} \not\thickapprox x^{n}),
\]
for each (equiv.\ some) $q \in (0, 1) \cap \mathbb{Q}$ and each $p \in [0, 1)\cap \mathbb{Q}$, $n \in \omega$ such that $\sqrt[n]{p}$ is irrational.
\end{Theorem}

Theorem \ref{Thm:base-product} provides a decision procedure for admissibility in $\rprlog$. 
It follows from the rational root theorem (see e.g. \cite[Thm.\ III.6.8]{Hu74}) that the set $\{ \langle n,p\rangle: n \in \omega\mbox{ and }p \in [0, 1]\cap\mathbb{Q}\mbox{ and }\sqrt[n]{p}\in\mathbb{Q}\}$ is decidable. Hence, the base for the admissible rules from Theorem \ref{Thm:base-product} forms a decidable set.
To determine whether a rule $\gamma_{1}, \dots, \gamma_{n} \rhd \varphi$ is admissible in $\rprlog$, 
we use the procedure consisting in enumerating all proofs  with the assumptions among $\gamma_{1}, \dots, \gamma_{n}$ in the logic obtained from $\rprlog$ by adding the  rules in Theorem \ref{Thm:base-product}  (and accepting if $\varphi$ is obtained)
and, simultaneously, enumerating all tuples $\vec{a}$ of elements in $[0, 1] \cap \mathbb{Q}$ such that $\gamma_{1}^{\rationalPQ}(\vec{a}) = \dots = \gamma_{n}^{\rationalPQ}(\vec{a})= 1$ (and rejecting if one is found such that $\varphi^{\rationalPQ}(\vec{a}) \ne 1$).

\begin{Corollary}\label{Cor:decid-RPLOG-admiss}
Admissible rules in $\rprlog$ form a decidable set.
\end{Corollary}

The next results (which will be proved later on) present a full characterization of  (hereditarily, actively, passively) 
structurally complete extensions of $\rprlog$. 
 
\begin{Theorem}\label{Thm:SC-product}
An extension $\vdash$ of $\rprlog$ is SC if and only if one of the following holds:
\benroman
\item\label{item:SC-product1} $\vdash$ is the structural completion of $\rprlog$ and, therefore,  it is algebraized by $\QQQ(\rationalPQ)$; or
\item\label{item:SC-product2} $\vdash$ is algebraized by one of the three proper subvarieties of $\class{RPA}$, in which case $\vdash$ is HSC.
\eroman
\end{Theorem}

\begin{Corollary}
\label{Cor:HSC-product}
An extension  of\/ $\rprlog$ is HSC if and only if it is SC.
\end{Corollary}

\begin{Corollary}
\label{Cor:ASC-product}
An extension of\/ $\rprlog$ is ASC if and only if it is SC.
\end{Corollary}

\begin{Corollary}\label{Cor:PSC-product}
An extension of\/ $\rprlog$ is PSC if and only if it is SC or it validates all
rules of the form
\[
\boldsymbol{c}_{q} \lor \bigvee_{i=1}^k (x_{i}^{n_i} \leftrightarrow \boldsymbol{c}_{p_{i}}) \rhd 0,
\]
where $k$ is a non-negative integer, $n_1,\ldots,n_k$ are natural numbers,  and  $q,p_1,\ldots,p_k \in [0, 1) \cap \mathbb{Q}$ are such that all numbers $\sqrt[n_i]{p_i}$ are irrational.
\end{Corollary}

\begin{Remark}
In view of Corollaries \ref{Cor:HSC-product} and \ref{Cor:ASC-product}, the notions of ASC, SC, and HSC are equivalent for extensions of $\rprlog$.  We will show that this equivalence cannot be extended to PSC, as there exist extensions of $\rprlog$ that are PSC, but not SC. 

%To this end,  observe that, for arbitrary nontrivial algebras $\A$ and $\B$ such that that $\QQQ(\A)$ is PSC and $\B\in\VVV(\A)$, the quasivariety $\QQQ(\A\times\B)$ is also PSC. To prove this,  notice that since $\QQQ(\A)$ and $\QQQ(\A \times \B)$ have the same free algebras, the two quasivarieties have the same passive quasiequations.  Then consider a quasiequation $\Phi$ that is passive in $\QQQ(\A \times \B)$. As $\Phi$ is also passive in $\QQQ(\A)$, its antecedent cannot be satisfied in the free algebras of $\QQQ(\A)$.  Since $\QQQ(\A)$ is PSC,  we can apply Theorem \ref{Thm:Bergman}(\ref{item:Ber-3}) obtaining that the antecedent of $\Phi$ cannot be satisfied in $\A$ and, therefore, in $\A \times \B$ either.  It follows that  $\A \times \B \vDash \Phi$ and, therefore,  $\QQQ(\A \times \B)\vDash \Phi$.  Hence, we conclude that $\QQQ(\A \times \B)$ is PSC.

%This observation allows us to find an extension of $\rprlog$ that is PSC but not SC. 

To this end,  consider the algebra $\rationalPQ\times \rationalPQ'$, where $\rationalPQ'$ is the expansion of the product algebra 
 $\rationalP$ in which all constants $\boldsymbol{c}_p$ with $p \ne 0$ are interpreted to the maximum element $1$ and in which $\boldsymbol{c}_0$ is interpreted as $0$.  Moreover, let $\vdash$ be the unique extension $\vdash$ of $\rprlog$ algebraized by $\QQQ(\rationalPQ\times \rationalPQ')$. Notice that $\vdash$  is PSC, by Corollary \ref{Cor:PSC-product}.  On the other hand,  $\QQQ(\rationalPQ\times \rationalPQ')$ is neither a proper subvariety of $\class{RPA}$ (because it contains $\rationalPQ$) nor $\QQQ(\rationalPQ)$ (because $\rationalPQ' \notin \QQQ(\rationalPQ)$). Therefore, $\vdash$ is not SC, by Theorem \ref{Thm:SC-product}.
 \qed
\end{Remark}

We begin by showing how to derive Theorem \ref{Thm:SC-product} from Theorem \ref{Thm:base-product}. 
To this end, suppose that Theorem \ref{Thm:base-product} holds. 
Recall that $\Fm_{\class{RPA}}(\omega))$ denotes  the free denumerably generated rational product algebra.

\begin{Lemma}
\label{Lem:: SC of RPA = Q(Q_P)}
The following holds for a subquasivariety  $\class{K}$ of $\class{RPA}$.
\benroman
\item{\label{item:sc_rpa_i}} $\rationalPQ\in\class{K}$ if and only if $\VVV(\class{K})=\class{RPA}$.
\item{\label{item:sc_rpa_ii}} If $\rationalPQ\in\class{K}$, then $\rationalPQ\leq \Fm_{\class{K}}(\omega)$.
\item{\label{item:sc_rpa_iii}} If $\rationalPQ\in\class{K}$, 
then $\QQQ(\rationalPQ) = \QQQ( \Fm_{\class{K}}(\omega))$.
\eroman
\end{Lemma}

\begin{proof}
(\ref{item:sc_rpa_i}): Suppose $\rationalPQ\in\class{K}$. Then $\VVV(\rationalPQ)\subseteq\VVV(\class{K})$. By Proposition \ref{Prop:Zuzka-trick}, $\class{RPA}=\VVV(\rationalPQ)$, therefore $\VVV(\class{K})=\class{RPA}$.
Now, let $\VVV(\class{K})=\class{RPA}$.  By Proposition \ref{Prop:Zuzka-trick},  the free algebras of $\class{RPA}$, $\class{K}$, and $\QQQ(\rationalPQ)$ coincide. 
% see e.g.  \cite[Thm. II.10.12]{BuSa00}  or  \cite[Theorem 14.1.2]{Ma73d}. 
In particular, since $\rationalPQ$ is the zero-generated free algebra in $\class{RPA}$, we conclude that $\rationalPQ\in\class{K}$.

(\ref{item:sc_rpa_ii}): This follows from the fact that if $\rationalPQ \in \class{K}$, then $\rationalPQ$ is the zero-generated free algebra of $\class{K}$.

(\ref{item:sc_rpa_iii}): By (\ref{item:sc_rpa_ii}), $\QQQ(\rationalPQ) \subseteq \QQQ( \Fm_{\class{K}}(\omega))$.
On the other hand, since the free algebras of $\QQQ(\rationalPQ)$ and of $\class{K}$ coincide, we have $ \Fm_{\class{K}}(\omega) \in \QQQ(\rationalPQ)$, hence 
$\QQQ( \Fm_{\class{K}}(\omega) ) \subseteq \QQQ(\rationalPQ)$.
\end{proof}

\begin{proof}[Proof of Theorem \ref{Thm:SC-product}]
Let $\class{K}_{\vdash}$ be the quasivariety of rational product algebras
algebraizing $\vdash$.

Suppose first that
 $\VVV(\class{K}_\vdash)=\class{RPA}$.
By Theorem \ref{Thm:Bergman}(\ref{item:Ber-1}), $\vdash$ is SC if and only if 
$\class{K}_\vdash=\QQQ(\Fm_{\class{K}_\vdash}(\omega))$.
Moreover,  applying  Lemma \ref{Lem:: SC of RPA = Q(Q_P)},   
 $\vdash$ is SC if and only if  $\class{K}_\vdash=\QQQ(\rationalPQ)$;
by Theorem \ref{Thm:base-product}, 
the last equality holds if and only if $\vdash$ is the structural completion of $\rprlog$.

Suppose on the other hand that $\VVV(\class{K}_\vdash)$ is a proper subvariety of $\class{RPA}$. By Theorem \ref{Thm:rational-variety}, this guarantees that $\class{K}_{\vdash}$ is term-equivalent to a quasivariety of product algebras. By theorem \ref{Thm:no-constants-product}, the variety of product algebras is primitive. It follows that $\class{K}_\vdash$ is a variety, i.e., the condition (\ref{item:SC-product2}) holds. In view of Theorem \ref{Thm:Bergman}(\ref{item:Ber-2}), $\vdash$ is HSC.
\end{proof}

\begin{proof}[Proof of Corollary \protect{\ref{Cor:HSC-product}}]
It is enough to show that all SC extensions of $\rprlog$ listed in Theorem \ref{Thm:SC-product} are HSC. For the extensions in item (\ref{item:SC-product2}), 
this fact follows from Theorems 
\ref{Thm:no-constants-product} and
\ref{Thm:rational-variety} (cf. the proof of Theorem \ref{Thm:SC-product}).

For the structural completion of $\rprlog$, we show that $\QQQ(\rationalPQ)$ is a minimal quasivariety.  To this end, observe that if $\A \in \QQQ(\rationalPQ)$ is nontrivial, then $\A$ validates the quasiequations of the form $\boldsymbol{c}_q \thickapprox 1 \Longrightarrow 0 \thickapprox 1$, for $q \in \Rat \cap[0, 1)$. Consequently, $\boldsymbol{c}_q^\A < 1$, for all $q\in[0,1)\cap\Rat$ and, therefore,  $\rationalPQ$ embeds into $\A$. This implies that $\QQQ(\rationalPQ) \subseteq \QQQ(\A)$. Hence, we conclude that $\QQQ(\rationalPQ)$ is a minimal quasivariety.
\end{proof}

\begin{proof}[Proof of Corollary \protect{\ref{Cor:ASC-product}}]
Let $\vdash$ be an ASC extension of $\rprlog$ and let $\class{K}_\vdash$ be the quasivariety algebraizing $\vdash$. If $\class{K}_\vdash\models \boldsymbol{c}_p\thickapprox 1$ for some $p\in [0,1)\cap\Rat$ then, by Theorem \ref{Thm:rational-variety}, $\class{K}_\vdash$ is contained in a proper subvariety of $\class{RPA}$. Thus, by Theorem \ref{Thm:SC-product}, the quasivariety $\class{K}_\vdash$ is primitive and $\vdash$ is SC.

Otherwise,   $\rationalPQ \in \class{K}_{\vdash}$.  This implies $\VVV(\class{K}_\vdash)=\class{RPA}$, by Lemma \ref{Lem:: SC of RPA = Q(Q_P)}(\ref{item:sc_rpa_i}).  Consequently,  $\class{RPA}$ and $\class{K}_\vdash$ have the same free algebras and, therefore, the same admissible quasiequations. It follows that $\rprlog$ and $\vdash$ have the same admissible rules.  Therefore, in order to prove that $\vdash$ is SC, it suffices to show that the rules that are admissible in $\rprlog$ are derivable in $\vdash$.  Clearly,  it will be enough to show that the rules in the base of the admissible rules for $\rprlog$ presented in Theorem
\ref{Thm:base-product} are derivable in $\vdash$. 

To this end, consider any rule $\varphi(x) \lor z \rhd z$ in this base.  Since $\emptyset\vdash \varphi(x)\lor 1$, this rule active in $\vdash$.  Furthermore, it is admissible, because $\rprlog$ and $\vdash$ have the same admissible rules. From the assumption  that $\vdash$ is ASC it follows that $\varphi(x)\lor z\vdash z$.  Hence, we conclude that $\vdash$ is SC, as desired.
\end{proof}

\begin{proof}[Proof of Corollary \ref{Cor:PSC-product}]
It suffices to show that a non SC extension $\vdash$ of $\rprlog$ is PSC if and only if it validates all
rules of the form
\begin{equation}\label{Eq:PSC-rules-new-writeup}
\boldsymbol{c}_{q} \lor \bigvee_{i=1}^k (x_{i}^{n_i} \leftrightarrow \boldsymbol{c}_{p_{i}}) \rhd 0,
\end{equation}
where $k$ is a non-negative integer, $n_1,\ldots,n_k$ are natural numbers,  and  $q,p_1,\ldots,p_k \in [0, 1) \cap \mathbb{Q}$ are such that all numbers $\sqrt[n_i]{p_i}$ are irrational.

Accordingly, consider a non SC extension $\vdash$ of $\rprlog$ and let $\class{K}_\vdash$ be the its equivalent algebraic semantics. By Theorems \ref{Thm:rational-variety} and \ref{Thm:SC-product}, we have $\VVV(\class{K}_\vdash)=\class{RPA}$; moreover by 
Lemma \ref{Lem:: SC of RPA = Q(Q_P)}(\ref{item:sc_rpa_i}), $\rationalPQ\in\class{K}_\vdash$. 
 
Let $T$ be the set of all formulas that are antecedents of one of the rules in (\ref{Eq:PSC-rules-new-writeup}). The stipulation that $\vdash$ validates the rules in (\ref{Eq:PSC-rules-new-writeup}) is equivalent to the following:
\benroman
\item\label{item:PSC proof cond i} Every nontrivial $\A \in \class{K}_\vdash$ validates 
all the sentences $\forall \vec x\; \delta(\vec x)\not\thickapprox 1$ with $\delta(\vec x)\in T$. 
\eroman

In view of Theorem \ref{Thm:Bergman}(\ref{item:Ber-3}),  the logic $\vdash$ is PSC if and only if the nontrivial members of $\class{K}_\vdash$ validate the same existential positive sentences. Since $\rationalPQ$ is free in $\class{K}_\vdash$, there exists a homomorphism from $\rationalPQ$ to any member of $\class{K}_\vdash$. Hence, every existential positive sentence which is valid in $\rationalPQ$ is also valid in $\class{K}_\vdash$. Consequently, the stipulation that $\vdash$ is PSC is equivalent to the following:
\benroman
\setcounter{enumi}{1}
\item\label{item:PSC proof cond ii} If $\A\in \class{K}_\vdash$ is nontrivial,  every existential positive sentence  valid in $\A$ is also valid in $\rationalPQ$.
\eroman
Therefore, it will be enough to show that conditions (\ref{item:PSC proof cond i}) and (\ref{item:PSC proof cond ii}) are equivalent.

Assume first that (\ref{item:PSC proof cond i}) holds and consider a nontrivial $\A \in \class{K}_\vdash$. Let 
\[
D=\{\delta^{\A}(\vec a) : \delta\in T \text{ and } \vec a \text{ is a tuple of elements in } A\}.
\]
By assumption, $1\not \in D$. We will verify that the set $D$ is closed under the join operation. Let  $\boldsymbol{c}_{q_1} \lor \epsilon_1(\vec x_1)$ and $\boldsymbol{c}_{q_2} \lor \epsilon_2(\vec x_2)$ be any terms in $T$ and $\vec a_1,\vec a_2$ tuples of elements in $A$. Put
$q \coloneqq \max({q_1,q_2})$. Then
\[
\boldsymbol{c}_{q_1} \lor \epsilon_1^\A(\vec a_1)\lor \boldsymbol{c}_{q_2} \lor \epsilon_1^\A(\vec a_2)=\boldsymbol{c}_{q} \lor \epsilon_1^\A(\vec a_1)\lor \epsilon_2^\A(\vec a_2)\in D.
\]  

By Lemma \ref{lem:: prime filter lemma}, there exists a prime filter $F$ of $\A$ such that $D\cap F=\emptyset$.  Furthermore, $\A/F$ is a chain,  by Theorem \ref{Thm: filters and congruences}. And since $D\cap F=\emptyset$, the algebra $\A/F$ validates all sentences $\forall \vec x\; \delta(\vec x)\not\thickapprox 1$ with $\delta\in T$. In particular, $\A/F$ validates all sentences listed in Theorem \ref{Thm:axiom-chain-Q} and, using this theorem, 
 $\A/F$ validates the universal theory of $\rationalPQ$. 

 To prove (\ref{item:PSC proof cond ii}), we will reason by contraposition. Accordingly,  consider a positive existential sentence $\Psi$ that fails in $\rationalPQ$.\ 
Then its negation $\boneg\Psi$ is equivalent to a universal 
sentence valid in $\rationalPQ$.\ Using the fact just proved, we infer that $\boneg\Psi$ is valid in $\A/F$. Consequently, $\Psi$ does not hold in $\A / F$. Since $\A / F$ is a homomorphic image of $\A$ and $\Psi$ is a positive existential sentence (and, therefore, it persists in homomorphic images), $\Psi$ is not valid in $\A$, thus establishing condition (\ref{item:PSC proof cond ii}).

Now assume (\ref{item:PSC proof cond ii}) holds. Let $\delta(\vec x)\in T$ and let $\A$ be a nontrivial algebra in $\class{K}_\vdash$. By the specification imposed on the parameters in $\delta$ and the fact that $\rationalPQ$ is a chain, the sentence $\exists \vec x\; \delta(\vec x)\thickapprox 1$
fails in $\rationalPQ$. By assumption, this implies that it also fails in $\A$. Hence, $\forall \vec x\; \delta(\vec x)\not\thickapprox 1$ holds in $\A$. This shows that the condition (\ref{item:PSC proof cond i}) holds.
\end{proof}

In order to prove Theorem \ref{Thm:axiom-chain-Q}, let us first collect a few relevant facts. 
By an \emph{$\ell$-group} we denote a lattice-ordered abelian group; 
an \emph{$o$-group} is a totally ordered $\ell$-group.
There exists a categorical equivalence between the class of product chains and the class of 
$o$-groups (see \cite{Cignoli-Torrens:ProductLogic} for an extension to 
a broader class of product algebras and $\ell$-groups).
Let us formulate a relevant part of this equivalence.

\begin{Proposition}[\protect{\cite[Thm.\ 4.1.8]{Ha98}, \cite[Thm.\ 2]{Hajek-Godo-Esteva:ProductLogic}}]
\label{Prop::o-groups}
Let $\A$ be a nontrivial product chain. 
Then there exists an $o$-group $\Lambda(\A)$  such that its negative cone $\{g \in \Lambda(A)  : g\leq^{\Lambda(\A)} 1 \}$ coincides with $A\smallsetminus\{0\}$ and for every $a,b\in A\smallsetminus\{0\}$ we have
\begin{gather*}
1^\A=1^{\Lambda(\A)},\quad a\leq ^\A b \text{ if and only if }a\leq ^{\Lambda(\A)} b,\quad  a\cdot^\A b=a\cdot^{\Lambda(\A)} b,\\ a\to^\A b=
\begin{cases}
1^{\Lambda(\A)} &\text { if } a\leq^\A b\\
b\cdot ^{\Lambda(\A)} (a^{-1})^{\Lambda(\A)} &\text { otherwise.} 
\end{cases}
\end{gather*}
\end{Proposition} 
Notice that $\Lambda(\A)$ is unique up to isomorphism. For $\Lambda(\standardP)$ and $\Lambda(\rationalP)$
 we take the (multiplicative) groups of positive reals and positive rationals.  

We will rely on the following results.

\begin{Theorem}[\protect{\cite[Thm. 5.3]{Savicky-CEGN:ProductTruthConstants}}]
\label{thm:Salvicky}
Let $\A$ be a rational product chain such that $\boldsymbol{c}_q<^\A 1$ for some (equiv.\ every) $q\in(0,1)\cap\Rat$. 
Then $\A$  partially embeds into $\standardPQ$.
\end{Theorem}

\begin{Theorem}\protect{\cite[Thm.\ II.1.6]{Hu74}}\label{Thm:fact_gen_subgroup_basis}
Let $m\in \omega$.
Let $\bm{F}$ be a free abelian group of rank $m$ and $\bm{G}$ a nontrivial subgroup of $\bm{F}$. 
Then $\bm{G}$ is a free abelian group. 
Moreover, there exists a basis $\{e_1,\dots, e_m\}$  of $\bm{F}$
and positive integers $k\leq m$  and $d_1 ,\dots, d_k$ such that 
$\{e_1^{d_1}, e_2^{d_2},\dots, e_k^{d_k}\}$ is a basis of $\bm G$.
\end{Theorem}

We are now ready to prove Theorem \ref{Thm:axiom-chain-Q}.
\begin{proof}
The sentences in the statement are valid in $\rationalPQ$. Thus it suffices to prove that every rational product algebra $\A$ validating them also validates the universal theory of $\rationalPQ$.   To show this,  it suffices to prove that every such rational product algebra $\A$ partially embeds into $\rationalPQ$.  Although we will use Theorem \ref{thm:Salvicky}, 
the result could not be obtained by partially embedding first $\A$ into $\standardPQ$ and then $\standardPQ$ into $\rationalPQ$, 
because $\standardPQ$ cannot be partially embedded into $\rationalPQ$ (for instance, 
$\exists x\,( x^2\thickapprox \boldsymbol{c}_{\nicefrac{1}{2}}) $ holds in $\standardPQ$, but 
fails in $\rationalPQ$).

Accordingly, consider a rational product algebra $\A$ validating the sentences in the statement. Clearly, $\A$ is a chain, as it validates $\forall xy ( x \leq y \boor y \leq x)$. 
Moreover, since $\boldsymbol{c}_{q} \not\thickapprox 1$ for $q \in [0, 1) \cap \mathbb{Q}$, 
we may assume that $\rationalPQ \leq \A$. 
Let $\B$ be a finite partial subalgebra of $\A$. We will find an embedding $h \colon \B \to \rationalPQ$.
To this end, we may assume without loss of generality that $B$ contains $\boldsymbol{c}_{q}^{\A}$ 
for some $q \in (0, 1) \cap \mathbb{Q}$ and moreover
 that $1^{\A}\in B$ and $0^{\A}\not\in B$. 
 Let $\A^-$ and $\rationalP$ be the product algebra reducts of $\A$ and of $\rationalPQ$. 
Let also $\bm\Lambda(\A^-)$ and  $\bm\Lambda(\rationalP)$ be the $o$-groups  associated with $\A^-$ and $\rationalP$ respectively as in Proposition \ref{Prop::o-groups}. 
 Clearly $\rationalP\leq \A^-$, whence we may assume that ${\bm\Lambda}(\rationalP)\leq{\bm\Lambda}(\A^-)$.
 Moreover, $0^{\A} \notin B$ gives $B \subseteq {\bm\Lambda}(\A^{-})$.

Consider the set $\mathcal{C}$ of constants. An expansion $\Lambda(\A)$ with the constants from $\mathcal{C}$
of the $o$-group $\Lambda(\A^-)$  is obtained by interpreting 
each $c\in \mathcal{C}$ in the expansion with the element that interprets $c$ in $\A$.
The expansions $\Lambda(\standardPQ)$ and $\Lambda(\rationalPQ)$ are defined analogously from $\Lambda(\standardP)$ and $\Lambda(\rationalP)$.

\begin{Claim}
\label{Claim:: Salvicky expanded}
 $\Lambda(\A)$ is partially embeddable into $\Lambda(\standardPQ)$.
\end{Claim}

\begin{proof}[Proof of the Claim]
Let $D$ be any finite set of elements in $\Lambda(\A)$. We may assume, without loss of generality, that $D$ is closed under reciprocals. 
By Theorem \ref{thm:Salvicky}, there exists an embedding $g\colon D\cap A\to (0,1]$ 
 of the   finite partial subalgebra   of $\A$ with universe  $D\cap A$
into $\standardPQ$. Let $f\colon D\to (0,\infty)$ be given by
$ f(d) = g(d)$  if $d\leq 1$, otherwise $f(d) = g(d^{-1})^{-1}$.
Then $f$ is an embedding of the   finite partial subalgebra    of $\Lambda(\A)$ with universe $D$
into $\Lambda(\standardPQ)$.  For instance, consider $c,d\in D$ 
such that $c\leq 1$, $d>1$ and $c\cdot d \in D$. Then $g(c)=g(c\cdot d\cdot d^{-1})=g(c\cdot d)\cdot g(d^{-1})$. 
Thus, $f(c\cdot d)=g(c\cdot d)=g(c)\cdot g(d^{-1})^{-1}=f(c)\cdot f(d)$.
\end{proof}

Let $\bm{G}_{B}$ be the $o$-subgroup of ${\bm \Lambda}(\A^{-})$ generated by $B$. (Then $\bm{G}_{B}$ is generated by $B$ also as a group.)  Let also $\bm{Q}_{B}$ be the $o$-subgroup of $\BoldLambda(\A^{-})$ whose universe is the intersection of $G_{B}$ with the universe of $\BoldLambda(\rationalP)$. Clearly, $\bm{Q}_{B}$ is also a $o$-subgroup of $\boldsymbol{G}_{B}$. As $B$ is finite, $\bm{G}_{B}$ is finitely generated. Furthermore, since the multiplication in ${\bm \Lambda}(\A^{-})$ preserves the strict order $<$, $\bm{G}_{B}$ is torsion free. Thus, the group reduct of $\bm{G}_{B}$ is free. 
Lastly, $\boldsymbol{c}_{q} \in B \cap \mathbb{Q} \subseteq {Q}_{B}$, whence $\boldsymbol{Q}_{B}$ is nontrivial.  Applying Theorem \ref{Thm:fact_gen_subgroup_basis} we obtain a basis $\{e_{1}, \dots, e_{m}\}$ for the group reduct of $\bm{G}_{B}$, and positive integers $k\leq m$ and $d_{1}, \dots, d_{k}$ such that  $\{ e_{1}^{d_{1}}, \dots, e_{k}^{d_{k}} \}$ is a basis for the group reduct of $\bm{Q}_{B}$. 
 Since $\boldsymbol{G}_{B}$ is totally ordered and multiplication is order preserving, we can also assume that $e_{1}, \dots, e_{m} \leq 1^{\A}$.
%
%Since $\bm{G}_{B}$ is an $o$-group, either $g<1$ or $g^{-1}<1$ for every $g\in G_B\smallsetminus\{1\}$. Therefore we  assume that $e_{1}, \dots, e_{m} < 1$ and, consequently, $e_{1}, \dots, e_{m}\in A\smallsetminus\{0,1\}$. 

\begin{Claim}
\label{Claim:: ei are rational}
We have $e_1,\ldots,e_k\in(0,1)\cap\Rat$.
\end{Claim}

\begin{proof}[Proof of the Claim]
Let $i\in\{1,\ldots,k\}$ and $q\coloneqq e_i^{d_i}$. 
Since $q$ belongs to the basis for $\bm Q_B$ and $Q_B\subseteq \Rat$,  we obtain that $q$ is rational. 
Moreover, as  $e_i^{d_i}=\boldsymbol{c}_q^\A$, the sentence $\forall x \, (x^{d_i}\not\thickapprox \boldsymbol{c}_q)$ fails in $\A$, so it is not among the axioms in Theorem \ref{Thm:axiom-chain-Q}.
It follows that $\sqrt[d_i]{q}$ is rational. 
As $\rationalPQ\leq\A$,  this yields $\sqrt[d_i]{q}\in A$. We obtain that $e_i^{d_i}=q=(\sqrt[d_i]{q})^{d_i}$ in $\A$. 
 As taking powers is injective in product algebras, this yields $e_i=\sqrt[d_i]{q}\in \Rat$. 
\end{proof}

Therefore,  we may assume that $d_{i} = 1$,  for all $1\leq i\leq k$. 
Summarizing the situation, $\{ e_{1}, \dots, e_{m} \}$ is a basis for $\boldsymbol{G}_{B}$ while
$\{ e_{1}, \dots, e_{k} \}$ a basis for $\boldsymbol{Q}_{B}$.
Therefore, every element $b\in B$ can be represented as 
\[
b=e_1^{l^b_1}\cdot\cdots\cdot e_m^{l^b_m}
\]
for unique integers $l^b_1,\ldots,l^b_m$. Let 
\[
 l\coloneqq \max\{|l^b_i| : i\in\{1,\ldots,m\}\text{ and } b\in B\},
 \]
  where $|l^b_i|$ is the absolute value of $l^b_i$.  Moreover let $C_0=\{e_1^{j_1}\cdot\cdots\cdot e_m^{j_m} : j_i\in{\mathbb Z}\text{ and } 
|j_i|\leq l \text{ for all }i\leq m\}$ and let $\C_0$ be the finite partial subalgebra  of $\Lambda(\A)$ 
with universe $C_0$.
We have $B\subseteq C_0\subseteq G_{\B}$.
Let $f\colon C_0\to (0,\infty)$ be an embedding of $\C_0$ into $\Lambda(\standardPQ)$ 
as in Claim \ref{Claim:: Salvicky expanded}. 
Notice that if $q\in C_0\cap Q_{\B}$, then  $f(q)=q$.  This is because $C_0$ is closed under reciprocals
and, since $\boldsymbol{c}_{q'}^{\Lambda(\A)}$ is defined for $q'\coloneqq \min(q,\nicefrac{1}{q})$, we have
\[
f(\boldsymbol{c}_{q'}^{\Lambda(\A)}) = \boldsymbol{c}_{q'}^{\Lambda(\standardPQ)}=q' \, \, \, \,  \text{ and } \, \,  \, \, 
f(\nicefrac{1}{q'}) = \nicefrac{1}{f(q')}.
\]

Let $C = C_0 \cap A$ and $\C$ be the partial subalgebra of $\A$ on $C$.  Since $\C$ extends $\B$,  to conclude the proof, it suffices to show that $\C$ embeds into $\rationalPQ$. 
Define 
\[
\Delta\coloneqq \min\left\{{\large\nicefrac{s}{r}} : r,s\in f[C]\text{ and }r<s\right\} \mbox{\ \ and\ \ }
d\coloneqq \sqrt[2\cdot lm]{\Delta}.
\]
Notice that $d>1$. Clearly, $e_1,\dots,e_m \in C$.

For every $i\in\{1,\ldots,m\}$ we pick a number $r_i$ in  $(0,1)\cap\Rat$ as follows. 
For $i\leq k$, let $r_i\coloneqq e_i=f(e_i)$. This is possible, because $e_i$ is rational by Claim \ref{Claim:: ei are rational}. For $i>k$, let $r_i$ be any element in $(0,1)\cap\Rat$ subject to the bounds $\nicefrac{r_i}{f(e_i)}<d$ and $\nicefrac{f(e_i)}{r_i}<d$.
It follows from the density of $\Rat$ in $\Real$ that such an $r_i$ exists.
   
Finally, we define a function $h\colon C\to (0,1]$ by putting, for $c\in C$,
\[
h(c)\coloneqq r_1^{l^c_1}\cdot\cdots\cdot r_m^{l ^c_m}.
\]    
We will verify that $h$ is an embedding of $\C$ into $\rationalPQ$.

\begin{Claim}
For every $q\in(0,1]\cap \Rat$ such that $\boldsymbol{c}^\A_q\in C$, we have $h(\boldsymbol{c}^\A_q)=q$.
\end{Claim}

\begin{proof}[Proof of the Claim]
Since $\boldsymbol{c}^\A_q \in C$ is rational and in $G_{\B}$, 
there are unique $l_1,\dots,l_k$ such that  $\boldsymbol{c}^\A_q = e_1^{l_1}\cdots e_k^{l_k}$ 
with $|l_i|\leq l$ for $i\leq k$.   Since $\rationalPQ \leq \A$ and $e_1, \dots, e_k \in [0, 1] \cap \Rat$, this implies that $e_1^{l_1}\cdots e_k^{l_k} = q$, where multiplication is computed in $\rationalPQ$. As $h$ is the identity map on the set $\{e_1, \dots, e_k\}$, the statement follows.
\end{proof}

\begin{Claim}
\label{Claim:order preservation}
The map $h$ is order preserving and injective on $C$. 
Consequently, $h[C]\subseteq (0,1]\cap\Rat$ and $h$ preserves the lattice operations. 
\end{Claim}

\begin{proof}[Proof of the Claim]
Consider an arbitrary $c\in C$.
Since all elements $e_1^{j_1}\cdots e_m^{j_m}$, where $|j_i|\leq l$, belong to $C_0$, 
by Claim \ref{Claim:: Salvicky expanded}
we have $f(c)=f(e_1)^{l_1^c}\cdots f(e_m)^{l_m^c}$. 
(However, these elements are not necessarily in $C$, which is why we extend 
partial embeddability to o-groups in Claim \ref{Claim:: Salvicky expanded}.) Thus,
\[
\frac{h(c)}{f(c)}=\left(\frac{r_1}{f(e_1)}\right)^{l^c_1}\cdot\cdots\cdot \left(\frac{r_m}{f(e_m)}\right)^{l ^c_m}< d^{lm}=\sqrt{\Delta}.
\] 
Similarly, 
\[
\frac{f(c)}{h(c)}<\sqrt{\Delta}.
\]
Applying the above inequalities and the definition of $\Delta$, 
for $c_1,c_2\in C$ such that $c_1<c_2$ we obtain
\[
h(c_1)<f(c_1)\cdot\sqrt{\Delta}
=\frac{f(c_1)\cdot\Delta}{\sqrt{\Delta}}\leq \frac{f(c_2)}{\sqrt{\Delta}}<h(c_2).
\]
Finally, as we assumed that $1^{\A}\in B\subseteq C$, it follows that $h(c)\leq h(1^{\A})=1$ for every $c\in C$. 
Since we assumed that $0^{\A}\not\in B$,  we conclude that $0 \not\in f[C]$. 
\end{proof}

\begin{Claim}
\label{Claim:multiplication preservation}
If $c_1,c_2,c_3\in C$ and $c_1\cdot^\A c_2=c_3$, then $h(c_1)\cdot^\rationalPQ h(c_2)=h(c_3)$.
\end{Claim}

\begin{proof}[Proof of the Claim.]
This follows from the uniqueness of the numbers $l^c_i$, where $c\in\{c_1,c_2,c_3\}$ and $i\in\{1,\ldots,m\}$.
\end{proof}

\begin{Claim}
If $c_1,c_2,c_3\in C$ and $c_1\to^\A c_2=c_3$, then $h(c_1)\to^\rationalPQ h(c_2)=h(c_3)$.
\end{Claim}

\noindent \textit{Proof of the Claim.}
If $c_3=1$, then $c_1\leq c_2$. Thus, by Claim \ref{Claim:order preservation}, 
$h(c_1)\leq h(c_2)$ and, therefore,  $h(c_1)\to^\rationalPQ h(c_2)=1=h(c_3)$. 
If $c_1> c_2$, by Proposition \ref{Prop::o-groups}, we have $c_3= c_2\cdot^{{\bm G}_B} (c_1^{-1})^{{\bm G}_B}$. Thus,  we have $l^{c_3}_i=l^{c_2}_i-l^{c_1}_i$ for $i\in\{1,\ldots,m\}$. 
Moreover, by Claim \ref{Claim:order preservation},  $h(c_1)> h(c_2)$. 
Therefore,  we obtain
\[
\pushQED{\qed} 
 h(c_3)= r_1^{l^{c_2}_1-l^{c_1}_1}\cdot\cdots\cdot\, r_m^{l^{c_2}_m-l^{c_1}_m}=
h(c_2)\cdot^{{\bm \Lambda}(\rationalPQ)} (h(c_1)^{-1})^{{\bm \Lambda}(\rationalPQ)}=
h(c_1)\to^{\rationalPQ} h(c_2).\qedhere \popQED
\]

Hence, $h \colon \C \to \rationalPQ$ is an embedding. This concludes the proof of \ref{Thm:axiom-chain-Q}.
\end{proof}

Lastly, we present a proof of Theorem \ref{Thm:base-product}.

\begin{proof}
Let $\vdash$ be the structural completion of $\rprlog$. 
By Theorem \ref{Thm:Bergman}(\ref{item:Ber-1}), $\vdash$ is the unique extension of $\rprlog$ 
algebraized by $\QQQ(\Fm_{\class{RPA}}(\omega))$, which 
by Lemma \ref{Lem:: SC of RPA = Q(Q_P)}(\ref{item:sc_rpa_iii}) equals  $\QQQ(\rationalPQ)$. Because of this, the problem of axiomatizing $\vdash$ relative to $\rprlog$ is equivalent to that of axiomatizing $\QQQ(\rationalPQ)$ relative to $\class{RPA}$. We will therefore focus on the latter. 

Let $T$ be the set of all terms $\boldsymbol{c}_p$ or $\boldsymbol{c}_q\leftrightarrow x^n$
with $p,q \in (0, 1)\cap\Rat$, $n\in\omega$, and $\sqrt[n]{q}$ irrational.  Moreover,  let $\A\in \class{RPA}$.
Clearly if $\A$ belongs to $\QQQ(\rationalPQ)$, then $\A$ validates all quasieqequations
 $\delta(x)\lor z\thickapprox 1\Longrightarrow z\thickapprox 1$ with $\delta(x)\in T$,
since they all hold in $\rationalPQ$ (recall that 
$\delta^{\rationalPQ} (a)<1$ for $\delta(x)\in T$,  $a\in[0,1]\cap\Rat$). 
Now assume, on the other hand, that $\A$ validates all these quasiequations.

\begin{Claim}
\label{Claim:: filters F_a}
Let $a\in A\smallsetminus{1}$. Then there exists a filter $F_a$ of $\A$ such that $a\not\in F_a$ and $\A/F_a$ validates the universal theory of $\rationalPQ$.
\end{Claim}

\begin{proof}[Proof of the Claim]
For $k \in \omega$, let 
\[
D_k\coloneqq\{a\lor \delta_1(b_1) \lor \cdots\lor \delta_k(b_k) :  \delta_1(x),\ldots,\delta_k(x)\in T\text{ and }b_1,\ldots b_k\in A\},
\]
and define $D\coloneqq\bigcup\{ D_k : k \in \omega\}$.

We will prove that $1 \notin D_k$,  by induction on $k$. As $a\neq 1$ and $D_0=\{a\}$,  we have $1\not\in D_0$. 
Suppose that $1\not\in D_{k-1}$ and, towards a contradiction, that $a\lor \delta_1(b_1) \lor \cdots\lor \delta_k(b_k)=1$ for some $\delta_1(x),\ldots,\delta_k(x)\in T$ and  $b_1,\ldots b_k\in A$. Since $\A$ satisfies the quasiequation $z\lor \delta_k(x)\thickapprox 1\Longrightarrow z\thickapprox 1$, we obtain that $a\lor \delta_1(b_1) \lor \cdots\lor \delta_{k-1}(b_{k-1})=1$ (consider an assignment which maps $z$ onto $a\lor \delta_1(b_1) \lor \cdots\lor \delta_{k-1}(b_{k-1})$ and $x$ onto $b_k$). This contradicts the assumption that $1\not \in D_{k-1}$. Hence, we conclude that
\[
1 \notin \bigcup_{k \in \omega}D_k = D.
\]

As $D$ is closed under the join operation and does not contain $1$, by Lemma \ref{lem:: prime filter lemma}, there exists a prime filter $F_a$ of $\A$ such that $F\cap D=\emptyset$. In particular, $a\not\in F_a$. It remains to show that $\A/F_a$ validates the sentences listed in Theorem  \ref{Thm:axiom-chain-Q}, i.e., that $\A/F_a$ is a chain and for every $\delta\in T$ and $b/F_a\in A/F_a$ we have $\delta(b/F_a)\neq 1$. By Theorem \ref{Thm: filters and congruences},  the primeness of $F_a$ yields that $\A / F_a$ is a chain. Then consider $\delta \in T$ and $b/F_a\in A/F_a$. By definition of $D$, we have $a \lor \delta(b) \in D$. As $F_a \cap D = \emptyset$ and $F_a$ is an upset,  this yields $\delta(b) \notin F_a$. Consequently, $\delta(b / F_a)  \ne 1$, as desired.
\end{proof}

For $a\in A\smallsetminus\{1\}$ let $F_a$ be a filter as in Claim \ref{Claim:: filters F_a}. Then all algebras $\A/F_a$ are in $\QQQ(\rationalPQ)$. Moreover, as $\bigcap\{F_a : a\in A\smallsetminus\{1\}\}=\{1\}$, the algebra $\A$ embeds into the product $\prod \{\A/F_a : a\in A\smallsetminus\{1\}\}$. This yields that $\A\in\QQQ(\rationalPQ)$.\footnote{The above argument can be replaced by the use of \cite[Cor.\ 3.8]{CzDz} 
or  \cite[Cor.~6]{Cintula-Noguera:Implicational2}. }
\end{proof}

\section{Extensions of rational G\"odel logic}\label{Sec:7}

The lattice of extensions of G\"odel logic $\mathbf{G}$ is notoriously transparent:

\begin{Theorem}[\protect{\cite{DzWr73}}]\label{Thm:GA-extensions}
Every extension of $\golog$ is axiomatic and the lattice of extensions of $\golog$ is a chain of order type $\omega + 1$. Consequently, $\golog$ is HSC.
\end{Theorem}
\noindent In algebraic parlance, the above result states that the quasivariety of G\"odel algebras is primitive, whence $\golog$ is HSC in view of Theorem \ref{Thm:Bergman}(\ref{item:Ber-2}). In this section, we shall see that the addition of rational constants to $\golog$ complicates the structure  of the lattice of (axiomatic) extensions of $\rgolog$.

\noindent  Given a real $r \in (0, 1]$, let $\boldsymbol{Q}_{r}$ be the rational G\"odel algebra with universe
\[
([0, r) \cap \mathbb{Q}) \cup \{ 1 \}
\]
The order relation of $\boldsymbol{Q}_{r}$ is the natural order in $\mathbb{Q}$. Accordingly, $\boldsymbol{Q}_{r}$ is a  chain. This settles the interpretation of the lattice connectives and of the implication (as for all $a, c \in Q_{r}$  we get $a \to c = 1$ if $a \leq c$, and $a \to c = c$ otherwise). Finally, given a rational $q \in [0, 1]$, the interpretation of of $\boldsymbol{c}_{q}$ in $\boldsymbol{Q}_{r}$ is $q$ if $q \in Q_{r}$, and $1$ otherwise. Notice that if $r=1$, then  $\boldsymbol{Q}_{1}=\boldsymbol{Q}_{G}$.

\noindent  Fix a denumerable set $\{ t_{n} : n \in \omega \}$ disjoint from $[0, 1]$.\ Given a rational $p \in [0, 1)\cap \mathbb{Q}$ and an ordinal $\gamma \in \omega +1$, let $\boldsymbol{Q}_{p}^{\gamma}$ be the rational G\"odel algebra with the universe
\[
([0, p] \cap \mathbb{Q}) \cup \{ 1 \} \cup \{ t_{n} : n < \gamma \}
\]
defined as follows. The order relation of $\boldsymbol{Q}_{p}^{\gamma}$ is given by the rule
\begin{align*}
a \leq c \Longleftrightarrow& \text{ either }c = 1 \text{ or }(a, c \in [0, 1]\text{ and }a \leq^{\mathbb{Q}} c) \text{ or }(a \in [0, 1) \text{ and }c \notin [0, 1])\\
& \text{ or }(a= t_{n} \text{ and }c = t_{m} \text{, for some }n \leq m).
\end{align*}
Accordingly, $\boldsymbol{Q}_{p}^{\gamma}$ is a chain. Similarly to the case of $\boldsymbol{Q}_{r}$, this settles the interpretation of the lattice connectives and of the implication. And given a rational $q \in [0, 1]$, the interpretation of $\boldsymbol{c}_{q}$ in $\boldsymbol{Q}_{p}^{\gamma}$ is $q$ if $q \in Q_{p}^{\gamma}$, and $1$ otherwise.

%It is easy to see that $\boldsymbol{Q}_{r}^{\gamma}$ is indeed a rational G\"odel algebra. Notice that $\boldsymbol{Q}_{1}^{0} = \rationalGQ$.

\begin{Theorem}\label{Thm:RGA-varieties} The following hold:
\benroman
\item\label{item:RGA-var-1} Every nontrivial variety $\class{K}$ of rational G\"odel algebras is of the form $\VVV(\boldsymbol{Q}_{r})$ for some $r \in (0, 1]$ or $\VVV(\boldsymbol{Q}_{p}^{\gamma})$ for some $\gamma \in \omega +1$ and $p \in [0, 1)\cap \mathbb{Q}$ . Furthermore, $\VVV(\boldsymbol{Q}_{r})$ is axiomatized by the equations $\{ \boldsymbol{c}_{q} \thickapprox 1 : q\in [r,1]\cap\mathbb{Q}\}$ and $\VVV(\boldsymbol{Q}_{p}^{\gamma})$ is axiomatized by the equations $\{ \boldsymbol{c}_{q} \thickapprox 1 : q \in (p, 1]\cap\mathbb{Q}\}$ and
\[
\Big(\bigvee_{0\leq i< j \leq n+2} (\boldsymbol{c}_{p} \lor x_{i}) \leftrightarrow(\boldsymbol{c}_{p} \lor x_{j})\Big) \thickapprox 1
\]
if  $\gamma = n \in \omega$, and by $\{ \boldsymbol{c}_{q} \thickapprox 1 : q \in (p, 1]\cap\mathbb{Q}\}$ otherwise.
\item\label{item:RGA-var-2} For all   $r_{1}, r_{2} \in (0, 1]$, $p_{1}, p_{2} \in [0, 1)\cap \mathbb{Q} $, and $\gamma_{1}, \gamma_{2} \in \omega +1$,
\begin{align*}
\VVV(\boldsymbol{Q}_{r_{1}}) \subseteq \VVV(\boldsymbol{Q}_{r_{2}}) &\Longleftrightarrow r_{1} \leq r_{2}, \\
\VVV(\boldsymbol{Q}_{r_{1}}) \subseteq \VVV(\boldsymbol{Q}_{p_{1}}^{\gamma_{1}}) &\Longleftrightarrow r_{1}\leq p_{1},\\
\VVV(\boldsymbol{Q}_{p_{1}}^{\gamma_{1}}) \subseteq \VVV(\boldsymbol{Q}_{r_{1}}) &\Longleftrightarrow p_{1}< r_{1},\\
\VVV(\boldsymbol{Q}_{p_{1}}^{\gamma_{1}}) \subseteq \VVV(\boldsymbol{Q}_{p_{2}}^{\gamma_{2}}) &\Longleftrightarrow\text{ either }p_{1} < p_{2} \text{ or }(p_{1} = p_{2} \text{ and }\gamma_{1} \leq \gamma_{2}).
\end{align*}

\item\label{item:RGA-var-3} $\mathcal{V}(\class{RGA})$ is an uncountable chain isomorphic to the poset obtained adding a new bottom element to the Dedekind--MacNeille completion of the lexicographic order of $[0, 1) \cap \mathbb{Q}$ and $\omega +1$.
\eroman
\end{Theorem}

\begin{Remark}
The axiomatization given in item (\ref{item:RGA-var-1}) can be simplified for varieties of the form $\VVV(\boldsymbol{Q}_{q})$ with $q \in \Rat \cap (0, 1]$, as these can be axiomatized by the single equation $\boldsymbol{c}_q
 \thickapprox 1$.  On the other hand, varieties of the form $\VVV(\boldsymbol{Q}_{r})$ with $r \in (0, 1] \smallsetminus \Rat$ do not admit a finite axiomatization.
\qed
\end{Remark}

 In view of the dual isomorphism between the lattice of axiomatic extensions of $\rgolog$ and $\mathcal{V}(\class{RGA})$, the above result provides a full description of the former as well.

Given a logic $\vdash$ and a set of formulas $\Sigma$, we denote by ${\vdash} + \Sigma$ the extension of $\vdash$ axiomatized relative to $\vdash$ by $\Sigma$.

\begin{Corollary}\label{Cor:axiomatic-extensions-RG}
Every consistent axiomatic extension of $\rgolog$ is of the form

\noindent $\rgolog_{r}\coloneqq\rgolog + \{c_{q} : q\in [r,1]\cap\mathbb{Q}\}$ for some $r\in(0,1]$,

 \noindent  $\rgolog_{p}^{\omega}\coloneqq\rgolog +\{c_{q} : q\in (p,1]\cap\mathbb{Q}\}$ for some rational $p\in[0,1)$ or

 \noindent    $\rgolog_{p}^{n}\coloneqq\rgolog_{p}^{\omega} + \bigvee_{0\leq i< j \leq n+2} (\boldsymbol{c}_{p} \lor x_{i}) \leftrightarrow(\boldsymbol{c}_{p} \lor x_{j})$ for some rational $p\in[0,1)$ and $n\in \omega$.

 \noindent   Moreover, the lattice of axiomatic extensions of $\rgolog$ is an uncountable chain dually isomorphic to the poset obtained adding a new bottom element to the Dedekind--MacNeille completion of the lexicographic order of $[0, 1) \cap \mathbb{Q}$ and $\omega +1$.
\end{Corollary}

On the other hand, the structure of the lattice of arbitrary extensions of $\rgolog$ is still largely unknown. For instance, the problem of determining whether the variety of rational G\"odel algebras is $\mathcal{Q}$-universal is still open. However, it is easy to see that it has uncountable chains and antichains. For chains, this is a consequence of Corollary \ref{Cor:axiomatic-extensions-RG}, while for antichains it suffices to notice that $\{ \QQQ(\boldsymbol{Q}_{r}) : r \in (0, 1] \}\cup\{ \QQQ(\boldsymbol{Q}_{p}^{0}) : p \in [0, 1)\cap \mathbb{Q} \}$ is a set of minimal quasivarieties (this can be proved by adapting the argument for the minimality of $\QQQ(\rationalPQ)$ in the proof of Corollary \ref{Cor:HSC-product}).

The rest of the section is dedicated to the proof of Theorem \ref{Thm:RGA-varieties}. We begin by the following observation:

\begin{Proposition}\label{Prop:RG-chains}
For every nontrivial rational G\"odel chain $\A$, there are $r \in (0, 1]$,  $p\in [0, 1)\cap \mathbb{Q }$, and $\gamma \in \omega + 1$ such that $\III\SSS\PPU(\A) = \III\SSS\PPU(\boldsymbol{Q}_{r})$ or  $\III\SSS\PPU(\A) = \III\SSS\PPU(\boldsymbol{Q}_{p}^{\gamma})$. Moreover, 
\benroman
\item\label{item:univ_god-alg1} $\III\SSS\PPU(\boldsymbol{Q}_{r})$ is axiomatized relative to the class of $\class{RGA}$ chains by the sentences
\[\boldsymbol{c}_{q'} \not\thickapprox 1 \mbox{ for all  } q'\in [0,r)\cap\mathbb{Q} \text{ and }   \boldsymbol{c}_{q} \thickapprox 1 \mbox { for all } q\in[r,1]\cap\mathbb{Q};\]
\item\label{item:univ_god-alg2}  $\III\SSS\PPU(\boldsymbol{Q}_{p}^{\omega})$ is  axiomatized relative to the class of $\class{RGA}$ chains by the sentences
\[
\boldsymbol{c}_{p} \not\thickapprox 1\text{ and }
 \boldsymbol{c}_{q} \thickapprox 1\mbox{ for all  } q\in(p,1]\cap\mathbb{Q};\]
\item\label{item:univ_god-alg3} $\III\SSS\PPU(\boldsymbol{Q}_{p}^{n})$ is  axiomatized relative to the class of $\class{RGA}$ chains by the sentences
\begin{align*}
\boldsymbol{c}_{p} \not\thickapprox 1, &\text{ }
 \boldsymbol{c}_{q} \thickapprox 1\mbox{ for all  } q\in(p,1]\cap\mathbb{Q},  \mbox{ and }\\
 \forall x_{0}\ldots x_{n+2}&\Big(\bigvee_{0\leq i< j \leq n+2} (\boldsymbol{c}_{p} \lor x_{i}) \leftrightarrow(\boldsymbol{c}_{p} \lor x_{j})\Big) \thickapprox 1.
\end{align*}
\eroman

\end{Proposition}

\noindent \textit{Proof.}
Let $\A$ be a rational G\"odel chain. We shall define an algebra $\boldsymbol{S_{A}}$ that embeds into $\A$. To this end, let $\C$ be the zero-generated subalgebra of $\A$. Since $\A$ is nontrivial, $\C = \boldsymbol{Q}_{r}$ for some $r \in (0, 1]$ or $\C = \boldsymbol{Q}_{p}^{0}$ for some $p\in [0, 1)\cap \mathbb{Q }$. If $\C = \boldsymbol{Q}_{r}$ for some $r \in (0, 1]$ then let $\boldsymbol{S_{A}} \coloneqq\C$. If $\C = \boldsymbol{Q}_{p}^{0}$ for some $p\in [0, 1)\cap \mathbb{Q }$, then let ${\downarrow}(C\smallsetminus \{ 1\})$ be the downset of $(C\smallsetminus \{ 1\})$ in $\A$. If $\omega \leq \vert A \smallsetminus {\downarrow}(C\smallsetminus \{ 1\}) \vert$, take $\boldsymbol{S_{A}}\coloneqq \boldsymbol{Q}_{p}^{\omega}$. While if $\vert A \smallsetminus {\downarrow}(C\smallsetminus \{ 1\}) \vert = n+1 \in \omega$, take $\boldsymbol{S_{A}} \coloneqq \boldsymbol{Q}_{p}^{n}$. In both cases, $\boldsymbol{S_{A}} \in \III\SSS(\A)$.   Therefore, in order to prove that $\boldsymbol{S_A}$ and $\A$ have the same universal theory, it suffices to show that $\A$ partially embeds into $\boldsymbol{S_{A}}$.

 To this end,  consider a finite partial subalgebra $\B$ of $\A$.  The elements of $B$ can be divided into those that are not the interpretation of any constant (denoted by $a_1, \dots, a_n$) and those that are (denoted by $\boldsymbol{c}_{q_{1}}^{\A}, \dots, \boldsymbol{c}_{q_{m}}^{\A} $). For the sake of simplicity,  we may assume that 
\[
0 = \boldsymbol{c}_{q_{1}}^{\A} < \boldsymbol{c}_{q_{2}}^{\A} < \dots < \boldsymbol{c}_{q_{m}}^\A = 1.
\]
As $\A$ is a chain, $[\boldsymbol{c}_{q_{1}}^{\A}, \boldsymbol{c}_{q_{2}}^{\A}), \dots, [\boldsymbol{c}_{q_{m-1}}^\A, \boldsymbol{c}_{q_{m}}^\A)$ is a partition of $A \smallsetminus \{ 1 \}$.  

Then consider the map $h \colon \B \to \boldsymbol{S_A}$ defined as follows.  For every $i \leq m-1$, let $a_{i_{1}} < \dots < a_{i_{k}}$ be the elements of $\{ a_{1}, \dots, a_{n} \}$ in the $i$-th component $[\boldsymbol{c}_{q_{i}}^{\A}, \boldsymbol{c}_{q_{i+1}}^{\A})$ of the above partition and choose some $b_{i_{1}}, \dots, b_{i_{k}} \in S_{A}$ 
such that
\[
\boldsymbol{c}_{q_{i}}^{\boldsymbol{S_{A}}} < b_{i_{1}} < \dots < b_{i_{k}} < \boldsymbol{c}_{q_{i+1}}^{\boldsymbol{S_{A}}}.
\]
If $i \ne m-1$, this is possible because $[\boldsymbol{c}_{q_{i}}^{\boldsymbol{S_{A}}}, \boldsymbol{c}_{q_{i+1}}^{\boldsymbol{S_{A}}})$ is an infinite set. 
While if $i = m-1$, this can be done by the construction of $\boldsymbol{S_{A}}$. Then 
let $h(a_{i_k}) \coloneqq b_{i_k}$. 
Furthermore, we set $h(\boldsymbol{c}_q^{\A}) = \boldsymbol{c}_{q}^{\A} = \boldsymbol{c}_{q}^{\boldsymbol{S_A}}$, 
for every $q \in \{ q_1, \dots, q_m \}$. This completes the definition of $h$. As in G\"odel chains the behaviour of the implication is fully determined by the order structure,  $h \colon \B \to \boldsymbol{S_A}$ is an embedding, as desired. We conclude that $\A$ and $\boldsymbol{S_A}$ have the same universal theory.

(\ref{item:univ_god-alg1}): Let $\boldsymbol{A}$ be a rational G\"{o}del chain validating the sentences in the statement. Then the zero-generated subalgebra  of $\A$ is $\boldsymbol{Q}_{r}$. Thus,  $\boldsymbol{S_{A}}=\boldsymbol{Q}_{r}$. Consequently, $\A$ and $\boldsymbol{Q}_{r}$ have the same universal theory and, in particular, $\A \in \III\SSS\PPU(\boldsymbol{Q}_{r})$.

(\ref{item:univ_god-alg2}): Let $\boldsymbol{A}$ be a rational G\"{o}del chain validating the sentences in the statement.  Then the zero-generated subalgebra of $\A$ is $\boldsymbol{Q}_{p}^{0}$. Thus,  $\boldsymbol{S_{A}}\in \{\boldsymbol{Q}_{p}^{\gamma}: \gamma\in\omega +1\}$.  Furthermore, since $\A$ and $\boldsymbol{S_A}$ have the same universal theory, $\A \in \III\SSS\PPU(\boldsymbol{S_A})$.  Thus,
\[
\boldsymbol{A}\in\III\SSS\PPU(\boldsymbol{S_{A}})\subseteq \III\SSS\PPU\SSS(\boldsymbol{Q}_{p}^{\omega}) =\III\SSS\PPU(\boldsymbol{Q}_{p}^{\omega}).
\]

(\ref{item:univ_god-alg3}): Let $\boldsymbol{A}$ be a rational G\"{o}del chain validating the sentences in the statement.  An argument similar to the one detailed for case (\ref{item:univ_god-alg2}) shows that $\boldsymbol{S_{A}}\in \{\boldsymbol{Q}_{p}^{\gamma}: \gamma\in\omega +1\}$.  Moreover, since $\A$ validates
\[
\forall x_{0}\ldots x_{n+2}\Big(\bigvee_{0\leq i< j \leq n+2} (\boldsymbol{c}_{p} \lor x_{i}) \leftrightarrow(\boldsymbol{c}_{p} \lor x_{j})\Big) \thickapprox 1,
\]
we obtain that $\vert A \smallsetminus {\downarrow}(C\smallsetminus \{ 1\}) \vert = m+1$ for some $m\leq n$, where $C$ is the universe of the zero-generated subalgebra of $\A$. Thus, $\boldsymbol{S_{A}}=\boldsymbol{Q}_{p}^{m}$, for some $m \leq n$. Consequently, 
\[
\pushQED{\qed}\boldsymbol{A}\in\III\SSS\PPU(\boldsymbol{S_{A}}) =\III\SSS\PPU(\boldsymbol{Q}_{p}^{m}) \subseteq \III\SSS\PPU\SSS(\boldsymbol{Q}_{p}^{n}) \subseteq\III\SSS\PPU(\boldsymbol{Q}_{p}^{n}).\qedhere \popQED
\]

\begin{Corollary}\label{Cor:RGA-generators}
Every variety of rational G\"odel algebras is generated by a set of algebras of the form $\boldsymbol{Q}_{r}$, where $r \in (0,1]$, 
 or  $\boldsymbol{Q}_{p}^{\gamma}$, where $p \in [0, 1)\cap \mathbb{Q}$ and $\gamma \in \omega + 1$.
\end{Corollary}

\begin{proof}
Recall that every variety is generated by its subdirectly irreducible members. As every subdirectly irreducible rational G\"odel algebra is a chain, the result follows from Proposition \ref{Prop:RG-chains}.
\end{proof}

\begin{proof}[Proof of Theorem \ref{Thm:RGA-varieties}.]

(\ref{item:RGA-var-2}): Consider   $r_{1}, r_{2} \in (0, 1]$, $p_{1}, p_{2} \in [0, 1)\cap \mathbb{Q} $ and $\gamma_{1}, \gamma_{2} \in \omega +1$. We need to prove that
\[\VVV(\boldsymbol{Q}_{r_{1}}) \subseteq \VVV(\boldsymbol{Q}_{r_{2}}) \Longleftrightarrow r_{1} \leq r_{2}.
\]

\noindent To prove the implication from left to right, we reason by contraposition. Accordingly, assume that  $r_{2} < r_{1}$. Since $\mathbb{Q}$ is dense in $\mathbb{R}$, there exists a rational $r_{2} \leq q < r_{1}$. Consequently, the equation $\boldsymbol{c}_{q} \thickapprox 1$ holds in $\boldsymbol{Q}_{r_{2}}$, but fails in $\boldsymbol{Q}_{r_{1}}$, whence $\boldsymbol{Q}_{r_{1}} \notin \VVV(\boldsymbol{Q}_{r_{2}})$. To prove the implication from right to left, if $r_{1} \leq r_{2}$,
then $\boldsymbol{Q}_{r_{1}} \in \HHH(\boldsymbol{Q}_{r_{2}}) \subseteq \VVV(\boldsymbol{Q}_{r_{2}})$.

Then we turn to prove that
\[
\VVV(\boldsymbol{Q}_{r_{1}}) \subseteq \VVV(\boldsymbol{Q}_{p_{1}}^{\gamma_{1}}) \Longleftrightarrow r_{1} \leq p_{1}.
\]
If $r_{1} \leq p_{1}$ then $\boldsymbol{Q}_{r_{1}} \in \HHH(\boldsymbol{Q}_{p_{1}}^{\gamma_{1}}) \subseteq \VVV(\boldsymbol{Q}_{p_{1}}^{\gamma_{1}})$.  If $p_{1} < r_{1}$, then there exists a rational $p_{1} < q < r_{1}$. Consequently, the equation $\boldsymbol{c}_{q} \thickapprox 1$ holds in $\boldsymbol{Q}_{p_{1}}^{\gamma_{1}}$, but fails in $\boldsymbol{Q}_{r_{1}}$, whence $\boldsymbol{Q}_{r_{1}} \notin \VVV(\boldsymbol{Q}_{p_{1}}^{\gamma_{1}})$.

Now, we will show that
\[\VVV(\boldsymbol{Q}_{p_{1}}^{\gamma_{1}}) \subseteq \VVV(\boldsymbol{Q}_{r_{1}}) \Longleftrightarrow p_{1} < r_{1}.
\]
 If $p_{1} < r_{1}$,  every finite partial subalgebra of $\boldsymbol{Q}_{p_{1}}^{\gamma_{1}}$ 
embeds  into some member of $\{\boldsymbol{Q}_{q} : q\in (p_{1},r_{1}]\cap \mathbb{Q}\}$. This implies that $\boldsymbol{Q}_{p_{1}}^{\gamma_{1}}$ validates the universal theory of $\{\boldsymbol{Q}_{q} : q\in (p_{1},r_{1}]\cap \mathbb{Q}\}$. Consequently, $\boldsymbol{Q}_{p_{1}}^{\gamma_{1}} \in \III\SSS\PPU(\{\boldsymbol{Q}_{q} : q\in (p_{1},r_{1}]\cap \mathbb{Q}\})$. As $\{\boldsymbol{Q}_{q} : q\in (p_{1},r_{1}]\cap \mathbb{Q}\} \subseteq \HHH(\boldsymbol{Q}_{r_{1}})$, this yields $\boldsymbol{Q}_{p_{1}}^{\gamma_{1}} \in \VVV(\boldsymbol{Q}_{r_{1}})$. 
If $r_{1} \leq p_{1}$, then $\boldsymbol{c}_{p_{1}} \thickapprox 1$ holds in $\boldsymbol{Q}_{r_{1}}$, but fails in $\boldsymbol{Q}_{p_{1}}^{\gamma_{1}}$ , whence $ \boldsymbol{Q}_{p_{1}}^{\gamma_{1}} \notin \VVV(\boldsymbol{Q}_{r_{1}})$.

Lastly, we will prove that 
\[
\VVV(\boldsymbol{Q}_{p_{1}}^{\gamma_{1}}) \subseteq \VVV(\boldsymbol{Q}_{p_{2}}^{\gamma_{2}}) \Longleftrightarrow\text{ either }p_{1} < p_{2} \text{ or }(p_{1} = p_{2} \text{ and }\gamma_{1} \leq \gamma_{2}).
\]
We prove the implication from left to right  by contraposition. Assume that either $p_{2} < p_{1}$ or ($p_{1} = p_{2}$ and $\gamma_{2} < \gamma_{1}$). First suppose that $p_{2} < p_{1}$, then there exists a rational $p_{2} < q < p_{1}$. Consequently, the equation $\boldsymbol{c}_{q} \thickapprox 1$ holds in $\boldsymbol{Q}_{p_{2}}^{\gamma_{2}}$, but fails in $\boldsymbol{Q}_{p_{1}}^{\gamma_{1}}$, whence $\boldsymbol{Q}_{p_{1}}^{\gamma_{1}} \notin \VVV(\boldsymbol{Q}_{p_{2}}^{\gamma_{2}})$. Then  suppose that $p_{1} = p_{2}=q\in[0,1)\cap \mathbb{Q}$ and $\gamma_{2} < \gamma_{1}$. Since $\gamma_{2} < \gamma_{1}$, necessarily $\gamma_{1} > 0$.  Furthermore, from $\gamma_{2} < \gamma_{1} \in \omega + 1$ it follows that $\gamma_{2} = n$ for some $n \in \omega$. Since $q < 1$ and $\gamma_{2} = n$, the interval $[\boldsymbol{c}_{q}, 1]$ in $\boldsymbol{Q}_{q}^{\gamma_{2}}$ is an $(n+2)$-element set, whence
\[
\boldsymbol{Q}_{p_{2}}^{\gamma_{2}} \vDash \bigvee_{0\leq i<j \leq n+2} (\boldsymbol{c}_{q} \lor x_{i}) \leftrightarrow(\boldsymbol{c}_{q} \lor x_{j}) \thickapprox 1.
\]
On the other hand, since $p_{1} = q$ and $\gamma_{1} > \gamma_{2} = n$, the interval $[\boldsymbol{c}_{q}, 1]$ in $\boldsymbol{Q}_{p_{1}}^{\gamma_{1}}$ has size $> n+2$. Consequently the above equation fails in $\boldsymbol{Q}_{p_{1}}^{\gamma_{1}}$, whence $\boldsymbol{Q}_{p_{1}}^{\gamma_{1}} \notin \VVV(\boldsymbol{Q}_{p_{2}}^{\gamma_{2}})$.

\noindent To prove the implication from right to left, if $p_{1} < p_{2}$, by previous items,  $ \VVV(\boldsymbol{Q}_{p_{1}}^{\gamma_{1}})\subseteq\VVV(\boldsymbol{Q}_{p_{2}})\subseteq \VVV(\boldsymbol{Q}_{p_{2}}^{\gamma_{2}})$.  If $p_{1} = p_{2}$ and $\gamma_{1} \leq \gamma_{2}$,
then $\boldsymbol{Q}_{p_{1}}^{\gamma_{1}} \in \SSS(\boldsymbol{Q}_{p_{2}}^{\gamma_{2}}) \subseteq \VVV(\boldsymbol{Q}_{p_{2}}^{\gamma_{2}})$.

(\ref{item:RGA-var-1}): Let $\class{K}$ be a nontrivial variety of rational G\"odel algebras. In view of Corollary \ref{Cor:RGA-generators}, $\class{K}$ is generated by a nonempty set of algebras $\{ \A_{i} \colon i \in I \}$ of the form $\boldsymbol{Q}_{r} $ or $\boldsymbol{Q}_{p}^{\gamma}$, where $r \in (0, 1]$, $p\in [0,1)\cap\mathbb{Q}$, and $\gamma \in \omega + 1$. We shall define an algebra $\boldsymbol{S}$ of the  previous type such that $\class{K}=\VVV(\boldsymbol{S})$.
Let
\[
s=\sup\left\{r \in [0, 1] : \boldsymbol{Q}_{r} \in\{ \A_{i} \colon i \in I \} \text{ or } \boldsymbol{Q}_{r}^{\gamma} \in\{ \A_{i} \colon i \in I \} \text{ for some } \gamma\in\omega + 1\right\}.
\]
If there exists $\gamma\in\omega +1$ such that $\boldsymbol{Q}_{s}^{\gamma} \in\{ \A_{i} \colon i \in I \}$,  let $\boldsymbol{S}\coloneqq \boldsymbol{Q}_{s}^{\delta}$ where $\delta= 	\sup\left\{\gamma\in\omega +1 : \boldsymbol{Q}_{s}^{\gamma} \in\{ \A_{i} \colon i \in I \}\right\}$.  Otherwise, let $\boldsymbol{S}\coloneqq \boldsymbol{Q}_{s}$.
By (\ref{item:RGA-var-2}), $\class{K}\subseteq \VVV(\boldsymbol{S})$. If $\boldsymbol{S}\in \{ \A_{i} \colon i \in I \}$, trivially $\VVV(\boldsymbol{S})\subseteq\class{K}$. If $\boldsymbol{S}\not\in \{ \A_{i} \colon i \in I \}$, then either $\boldsymbol{S}=\boldsymbol{Q}_{s}$ or $\boldsymbol{Q}_{s}^{\omega}$.     In both cases, every finite partial subalgebra in $\boldsymbol{S}$ embeds into some member of $\{ \A_{i} \colon i \in I \}$. As a consequence,  $\boldsymbol{S}$ validates the universal theory of $\{ \A_{i} \colon i \in I \}$, whence $\boldsymbol{S}\in\III\SSS\PPU(\{ \A_{i} \colon i \in I \})\subseteq\class{K}$. We conclude that $\class{K} = \VVV(\boldsymbol{S})$ and, therefore, that every variety of rational G\"odel algebras is generated by an algebra of the form $\boldsymbol{Q}_r$ or $\boldsymbol{Q}_p^\gamma$.

In order to axiomatize varieties of the form  $\VVV(\boldsymbol{Q}_{p}^{\gamma})$, let $\Sigma$ be the set of equations given by the statement. First observe that $\boldsymbol{Q}_{p}^{\gamma} \vDash \Sigma$. Then consider a rational G\"odel algebra $\A \notin \VVV(\boldsymbol{Q}_{p}^{\gamma})$. As we showed in the above paragraph, the variety $\VVV(\A)$ is generated by an algebra of the form of the form $\boldsymbol{Q}_{r}$ for some $r\in(0,1]$ or  $\boldsymbol{Q}_{p'}^{\delta}$ for some $p' \in [0, 1)\cap \mathbb{Q}$ and $\delta \in \omega + 1$. If $\VVV(\A)=\VVV(\boldsymbol{Q}_{r})$, since $\A \notin \VVV(\boldsymbol{Q}_{p}^{\gamma})$,  we get $\VVV(\boldsymbol{Q}_{r}) \nsubseteq \VVV(\boldsymbol{Q}_{p}^{\gamma})$. By (\ref{item:RGA-var-2}),  $p<r$.  Thus, there is a rational $p<q<r$ such that $c_{q}^{\boldsymbol{Q}_{r}}\neq 1$.  In that case, $\boldsymbol{Q}_{r} \nvDash \Sigma$, whence $\A \nvDash \Sigma$. If $\VVV(\A)=\VVV(\boldsymbol{Q}_{p'}^{\delta})$, since $\A \notin \VVV(\boldsymbol{Q}_{p}^{\gamma})$,  we get $\VVV(\boldsymbol{Q}_{p'}^{\delta}) \nsubseteq \VVV(\boldsymbol{Q}_{p}^{\gamma})$. By (\ref{item:RGA-var-2}), either $p<p'$ or $p=p'$ and $\gamma<\delta$. If $p<p'$, similar to the previous case $\boldsymbol{Q}_{p'}^{\delta} \nvDash c_{p'}\approx 1$, whence $\A \nvDash \Sigma$. If $p=p'$ and $\gamma<\delta$, then $\gamma=n\in\omega$ and
\[
\boldsymbol{Q}_{p'}^{\delta} \nvDash \Big(\bigvee_{0\leq i, j \leq n+2} (\boldsymbol{c}_{p} \lor x_{i}) \leftrightarrow(\boldsymbol{c}_{p} \lor x_{j})\Big) \thickapprox 1,
\]
whence $\A \nvDash \Sigma$. Thus, we conclude that $\Sigma$ axiomatizes $\VVV(\boldsymbol{Q}_{p}^{\gamma})$.

 It only remains to axiomatize varieties of the form $\VVV(\boldsymbol{Q}_{r})$ for $r\in(0,1]$.   Since $c_{q}^{\boldsymbol{Q}_{q}}=1$ for every rational $q\in[r,1]$,  the equations in the statement are valid in $\VVV(\boldsymbol{Q}_{r})$. A similar argument as in the case of  varieties of the form $\VVV(\boldsymbol{Q}_{p}^{\gamma})$ shows that if $\A\not\in\VVV(\boldsymbol{Q}_{r})$, then there is a rational $r\leq q<1$ such that $\A \nvDash \boldsymbol{c}_{q}\thickapprox 1$.

(\ref{item:RGA-var-3}): Let $\mathcal{V}(\class{RGA})^{-}$ be the poset of nontrivial varieties of rational G\"odel algebras. Notice that $\mathcal{V}(\class{RGA})^{-}$ is indeed a complete lattice and that $\mathcal{V}(\class{RGA})$ is obtained adding a new bottom element to $\mathcal{V}(\class{RGA})^{-}$. Therefore, to conclude the proof, it suffices to show that $\mathcal{V}(\class{RGA})^{-}$ is isomorphic to the Dedekind--MacNeille completion \cite{McNe37} of the poset $\mathbb{X}$ obtained by endowing the direct product
\[
([0, 1) \cap \mathbb{Q}) \times (\omega +1)
\]
with the lexicographic order of $\mathbb{Q}$ and $\omega +1$. By (\ref{item:RGA-var-1}) and (\ref{item:RGA-var-2}), the map $f \colon \mathbb{X} \to \mathcal{V}(\class{RGA})^{-}$, defined by   
\[
f(\langle q, \gamma \rangle) \coloneqq
\left\{
  \begin{array}{ll}
    \VVV(\boldsymbol{Q}_{0}^{\gamma}), & \hbox{if }q=0; \\
    \VVV(\boldsymbol{Q}_{q}) & \hbox{if }q\neq 0 \hbox{ and } \gamma=0; \\
   \VVV(\boldsymbol{Q}_{q}^{n}) & \hbox{if } q\neq 0 \hbox{ and } \gamma=n+1;\\
    \VVV(\boldsymbol{Q}_{q}^{\omega}) & \hbox{if } q\neq 0 \hbox{ and } \gamma=\omega,
  \end{array}
\right.
\]
is an order embedding. Furthermore, $f[\mathbb{X}]$ is both join-dense and meet-dense in the complete lattice $\mathcal{V}(\class{RGA})^{-}$. As, up to isomorphism, the Dedekind--MacNeille completion of a poset $\mathbb{Y}$  is the only completion in which $\mathbb{Y}$ is both join-dense and meet-dense \cite[Prop.\ 1]{BanBrun67} (see also \cite{Bruns62}), we conclude that $\mathcal{V}(\class{RGA})^{-}$ is isomorphic to the Dedekind--MacNeille completion of $\mathbb{X}$, as desired.
\end{proof}

\section{Structural completeness in rational G\"odel logic}\label{Sec:8}

It is well known that $\golog$ is HSC \cite{DzWr73}.\ While this is false for $\rgolog$, it is still possible to obtain a full characterization of structural completeness and its variants in extensions of $\rgolog$. The next result characterizes PSC extensions of $\rgolog$.

\begin{Theorem}\label{Thm:RGA-PSC}
The following are equivalent for an extension $\vdash$  of $\rgolog$:
\benroman
\item\label{RGA-PSC1} $\vdash$ is PSC;
\item\label{RGA-PSC2} $\vdash$ is algebraized by a quasivariety with the JEP;
\item\label{RGA-PSC3} $\vdash$ is algebraized by a quasivariety whose nontrivial members have isomorphic zero-generated subalgebras.
%\item\label{RGA-PSC3} $\vdash$ is algebraized by a quasivariety $\class{K}$ such that every two nontrivial members of $\class{K}$ have isomorphic zero-generated subalgebras.
\eroman
\end{Theorem}

The other variants of structural completeness turn out to be equivalent among extensions of $\rgolog$.

\begin{Theorem}\label{Thm:RGA-SC}
The following are equivalent for an extension $\vdash$  of $\rgolog$:
\benroman
\item\label{RGA-SC1} $\vdash$ is HSC;
\item\label{RGA-SC2} $\vdash$ is SC;
\item\label{RGA-SC3} $\vdash$ is ASC;
\item\label{RGA-SC4} $\vdash$ is algebraized by a quasivariety $\class{K}$ generated by a chain $\A$.

\eroman
Furthermore, in condition (\ref{RGA-SC4}) $\A$ can be chosen either trivial or of  the form $\boldsymbol{Q}_{r} $ or $\boldsymbol{Q}_{p}^{\gamma}$, where $r \in (0, 1]$, $p\in [0,1)\cap\mathbb{Q}$ and $\gamma \in \omega + 1$.
\end{Theorem}

\begin{proof}[Proof of Theorem \ref{Thm:RGA-PSC}.]
The implication (\ref{RGA-PSC1})$\Rightarrow$(\ref{RGA-PSC2}) is a consequence of Proposition \ref{Prop:PSC-JEP}, while (\ref{RGA-PSC2})$\Rightarrow$(\ref{RGA-PSC3}) is straightforward.

(\ref{RGA-PSC3})$\Rightarrow$(\ref{RGA-PSC1}): Let $\class{K}$ be the quasivariety algebraizing $\vdash$. In view of Theorem \ref{Thm:Bergman}(\ref{item:Ber-3}), it suffices to show that every two nontrivial members validate the same positive existential sentences. To this end, consider two nontrivial $\A, \B \in \class{K}$. By assumption the zero-generated algebra $\C$ of $\A$ and of $\B$ coincide.

By Lemma \ref{lem:: prime filter lemma}, the set ${\mathcal F}$ of prime filters $F$ of $\A$ such that $F \cap C = \{ 1 \}$ is nonempty,
and we can apply Zorn's lemma, obtaining a maximal $F \in \mathcal{F}$. 
% Notice that the filter $F$ is prime: suppose the contrary, with a view to contradiction. 
% Then there are $a, b \in A$ such that $a \lor b \in F$ and $a, b \notin F$. 
% Because of the maximality of $F$, there are $\boldsymbol{c}_p, 
% \boldsymbol{c}_q \in C \smallsetminus \{ 1 \}$ and $g_{a}, g_{b} \in F$ such that 
% $g_{a} \land a \leq \boldsymbol{c}_{p}$ and $g_{b} \land b \leq \boldsymbol{c}_{q}$. 
% By symmetry, we can assume that $p \leq q$, whence $\boldsymbol{c}_{p} \leq \boldsymbol{c}_{q}$. Consequently,
% \[
% (g_{a} \land g_{b}) \land (a \lor b) = (g_{a} \land g_{b} \land a) \lor (g_{a} \land g_{b} \land b) \leq \boldsymbol{c}_{q}.
% \]
% As $g_{a}, g_{b}, a \lor b \in F$, we get $(g_{a} \land g_{b}) \land (a \lor b) \in F$. 
% Together with the above display, this yields $\boldsymbol{c}_{q} \in F$, a contradiction. Thus, we conclude that $F$ is prime.
% 
% Then, consider the relation
% \[
% \theta \coloneqq \{ \langle a, c \rangle \in A \times A : a \to c, c \to a \in F \}.
% \]
% Since $F$ is a lattice filter, $\theta$ is a congruence on $\A$. 
As  $F$ is prime, by Theorem \ref{Thm: filters and congruences} $\A / F$ is a chain.
By construction of $F$ we know that the zero-generated subalgebra of $\A/ F$ is isomorphic to $\C$. Therefore, we may assume, without loss of generality, that $\C \leq \A / F$. 
Clearly, either $\boldsymbol{Q}_{r}$ or $\boldsymbol{Q}_{p}^{0}$ is the zero-generated subalgebra of $\A / F$, namely $\C$. 
Furthermore, by the maximality of $F$, we get that if $a \in A / F$ is strictly larger than 
all the elements of $C \smallsetminus \{ 1 \}$, then $a = 1$. 
Together with Proposition \ref{Prop:RG-chains}, this yields that $\III\SSS\PPU(\A/F ) = \III\SSS\PPU(\boldsymbol{Q}_{r})$ for some $r \in (0, 1]$ or  $\III\SSS\PPU(\A/F ) = \III\SSS\PPU(\boldsymbol{Q}_{p}^{0})$ for some $p \in [0, 1)\cap\mathbb{Q}$. 
Thus, $\A / F \in \III\SSS\PPU(\C)$. Since $\C \leq \B$, this yields $\A / F \in \III\SSS\PPU(\B)$. 
Then there exists an embedding $f_{A} \colon \A /F \to \B_{u}$, where $\B_{u}$ is an ultrapower of $\B$. 
Let $g_{A} \colon \A \to \A / F$ be the canonical surjection. 
Then the composition $h_{A} \colon f_{A} \circ g_{A}$ is a homomorphism from $\A$ to $\B_{u}$. Since positive existential sentences persist in homomorphic images,  extensions,  and ultraroots, we conclude that every positive existential sentence that is true of $\A$ is also true of $\B$.
\end{proof}

In order to prove Theorem \ref{Thm:RGA-SC}, we rely on the following observation.

\begin{Proposition}\label{Prop:RGA-embedding}
Let $\class{K}$ be a nontrivial subquasivariety of $\class{RGA}$ and $\Fm_{\class{K}}(\omega)$ its denumerably generated free algebra. Then there are $r \in [0, 1)$,  $p\in [0, 1)\cap \mathbb{Q }$ and $\gamma \in \omega + 1$ such that $\QQQ(\Fm_{\class{K}}(\omega)) = \QQQ(\boldsymbol{Q}_{r})$ or $\QQQ(\Fm_{\class{K}}(\omega)) = \QQQ(\boldsymbol{Q}_{p}^{\gamma})$.
\end{Proposition}

\noindent \textit{Proof.} By Theorem \ref{Thm:RGA-varieties}, there are $r \in [0, 1)$,  $p\in [0, 1)\cap \mathbb{Q }$, and $\gamma \in \omega + 1$ such that $\Fm_{\class{K}}(\omega)$ is the denumerably generated free algebra of $\VVV(\boldsymbol{Q}_{r})$ or $\VVV(\boldsymbol{Q}_{p}^{\gamma})$. This yields that $\Fm_{\class{K}}(\omega)$ is also the denumerably generated free algebra of $\QQQ(\boldsymbol{Q}_{r})$ or $\QQQ(\boldsymbol{Q}_{p}^{\gamma})$, whence $\QQQ(\Fm_{\class{K}}(\omega)) \subseteq \QQQ(\boldsymbol{Q}_{r})$ or $\QQQ(\Fm_{\class{K}}(\omega)) \subseteq \QQQ(\boldsymbol{Q}_{p}^{\gamma})$. If $\Fm_{\class{K}}(\omega)=\Fm_{\VVV(\boldsymbol{Q}_{r})}(\omega)$, then  $\boldsymbol{Q}_{r}$ is the zero-generated subalgebra of $\Fm_{\class{K}}(\omega)$. Whence $\QQQ(\boldsymbol{Q}_{r}) \subseteq \QQQ(\Fm_{\class{K}}(\omega))$.\ Similarly, if  $\Fm_{\class{K}}(\omega)=\Fm_{\VVV(\boldsymbol{Q}_{p}^{0})}(\omega)$, then  $\QQQ(\boldsymbol{Q}_{p}^{0}) \subseteq \QQQ(\Fm_{\class{K}}(\omega))$. Finally, assume $\Fm_{\class{K}}(\omega)=\Fm_{\VVV(\boldsymbol{Q}_{p}^{\gamma})}(\omega)$ with $\gamma > 0$. Notice that
\[
\boldsymbol{Q}_{p}^{\gamma} \in \III\SSS\PPU(\{ \boldsymbol{Q}_{p}^{n} : n \in \omega \text{ and }1 \leq n \leq \gamma \}).
\]
Consequently, to conclude the proof, it suffices to show that each $\boldsymbol{Q}_{p}^{n}$ (where $n \in \omega$ and $1 \leq n \leq \gamma)$ embeds into $\Fm_{\class{K}}(\omega)$. This can be done by a straightforward adaptation of the method described in \cite{DzWr73} for the case of G\"odel algebras without constants.

We shall sketch it for the sake of completeness. For every $1 \leq n \leq \gamma$ such that $n \in \omega$, the algebra $\boldsymbol{Q}_{p}^{n}$ is the chain consisting of the interval $[0, p] \cap \mathbb{Q}$ on top of which we added the $n+1$ element chain
\[
t_{0} < t_{1} < \dots < t_{n-1} < 1.
\]
Then $\boldsymbol{Q}_{p}^{n}$ can be embedded into $\boldsymbol{Fm}_{\class{K}}(\omega)$ using the map that is the identity on $([0, p]\cap \mathbb{Q}) \cup \{ 1 \}$ and that sends $t_{i}$ to the (equivalence class of) the formula $\varphi_{i}$, where
\[
\pushQED{\qed} \varphi_{0} \coloneqq \boldsymbol{c}_{p} \lor  x_{2} \lor (x_{2} \to x_{1}) \text{ and }\varphi_{j+1} \coloneqq x_{j+2} \lor (x_{j+2} \to \varphi_{j}).\qedhere \popQED
\]

\begin{proof}[Proof of Theorem \ref{Thm:RGA-SC}.]
The implications (\ref{RGA-SC1})$\Rightarrow$(\ref{RGA-SC2}) and (\ref{RGA-SC2})$\Rightarrow$(\ref{RGA-SC3}) are straightforward.

(\ref{RGA-SC3})$\Rightarrow$(\ref{RGA-SC2}): Let $\class{K}$ be the quasivariety algebraizing $\vdash$ and $\Fm_{\class{K}}(\omega)$ its denumerably generated free algebra. Suppose, with a view to contradiction, that $\vdash$ is not SC. As $\vdash$ is ASC, this means that it is not PSC. In view of Theorem \ref{Thm:RGA-PSC}, there is a zero-generated algebra $\C \in \class{K}$ different from the zero-generated subalgebra $\Fm_{\class{K}}(0)$ of $\Fm_{\class{K}}(\omega)$. Clearly, there is $q \in [0, 1) \cap \mathbb{Q}$ such that the equation $\boldsymbol{c}_{q} \thickapprox 1$ holds in $\C$, but not in $\Fm_{\class{K}}(\omega)$. 
Furthermore, by Proposition \ref{Prop:RGA-embedding} there are $r \in (0, 1]$, $p\in[0,1)\cap\mathbb{Q}$, and $\gamma \in \omega +1$ such that $\QQQ(\Fm_{\class{K}}(\omega)) = \QQQ(\boldsymbol{Q}_{r})$ or $\QQQ(\Fm_{\class{K}}(\omega)) = \QQQ(\boldsymbol{Q}_{p}^{\gamma})$. 
Let then $\A \in \{ \boldsymbol{Q}_{r}, \boldsymbol{Q}_{p}^{\gamma} \}$ be such that 
 $\QQQ(\Fm_{\class{K}}(\omega)) = \QQQ(\A)$.
In particular, $\A \nvDash \boldsymbol{c}_{q} \thickapprox 1$. 
As $\A$ is a chain, we have
\[
\A \vDash x \lor \boldsymbol{c}_{q} \thickapprox 1 \boto x \thickapprox 1.
\]
Since $\vdash$ is ASC, by Theorem \ref{Thm:Bergman}(\ref{item:Ber-4}) 
$\C \times \Fm_{\class{K}}(0) \in \QQQ(\Fm_{\class{K}}(\omega)) = \QQQ(\A)$, whence
\[
\C \times \Fm_{\class{K}}(0) \vDash x \lor \boldsymbol{c}_{q} \thickapprox 1 \boto x \thickapprox 1.
\]
But this is false, as witnessed by the assignment $x \longmapsto \langle 0, 1 \rangle$. Hence, we conclude that $\vdash$ is PSC and, therefore, SC.

(\ref{RGA-SC2})$\Rightarrow$(\ref{RGA-SC4}): Suppose that $\vdash$ is SC. By Theorem \ref{Thm:Bergman}(\ref{item:Ber-1}) and Proposition \ref{Prop:RGA-embedding}, $\vdash$ is algebraized by a quasivariety that is either trivial or of the form $\QQQ(\boldsymbol{Q}_{r})$ for some $r\in(0,1]$ or $\QQQ(\boldsymbol{Q}_{p}^{\gamma})$ for some $p\in[0,1)\cap\mathbb{Q}$ and $\gamma \in \omega +1$.

(\ref{RGA-SC4})$\Rightarrow$(\ref{RGA-PSC1}): If $\vdash$ is algebraized by the trivial quasivariety, then $\vdash$ is clearly HSC. Then we consider the case where it is algebraized by a quasivariety $\QQQ(\A)$ where $\A$ is a nontrivial chain. By Proposition \ref{Prop:RG-chains} we can assume that $\A = \boldsymbol{Q}_{r}$ for some $r \in (0, 1]$ or $\A = \boldsymbol{Q}_{p}^{\gamma}$ for some $p\in[0,1)\cap\mathbb{Q}$ and $\gamma \in \omega +1 $. 

Suppose first that $\A = \boldsymbol{Q}_{r}$. Since the zero-generated subalgebra of every nontrivial member of $\QQQ(\boldsymbol{Q}_{r})$ is $\boldsymbol{Q}_{r}$, the quasivariety $\QQQ(\boldsymbol{Q}_{r})$ is minimal and, therefore, $\vdash$ is HSC. 

Then we consider the case where $\A = \boldsymbol{Q}_{p}^{\gamma}$. In view of Theorem \ref{Thm:Bergman}(\ref{item:Ber-1}), to prove that $\vdash$ is HSC, it suffices to show that every subquasivariety of $\QQQ(\boldsymbol{Q}_{p}^{\gamma})$ is generated as a quasivariety by its denumerably generated free algebra. This is true for $\QQQ(\boldsymbol{Q}_{p}^{\gamma})$, by Proposition \ref{Prop:RGA-embedding}. Then consider a proper subquasivariety $\class{K}$ of $\QQQ(\boldsymbol{Q}_{p}^{\gamma})$. Since $\class{K}$ is nontrivial,  $\boldsymbol{Q}_{p}^{0} \in \class{K}$. On the other hand, as $\class{K}$ is proper, $\boldsymbol{Q}_{p}^{\gamma} \notin \class{K}$. Consequently, there is $n \in \omega$ such that for all $m \in \omega +1 $,
\[
\boldsymbol{Q}_{p}^{m} \in \class{K} \Longleftrightarrow m \leq n.
\]
Together with Proposition \ref{Prop:RGA-embedding}, this implies $\QQQ(\Fm_{\class{K}}(\omega)) = \QQQ(\boldsymbol{Q}_{p}^{n})$. In particular,
\[
\class{K}\vDash \Big(\bigvee_{0\leq i< j \leq n+2} (\boldsymbol{c}_{p} \lor x_{i}) \leftrightarrow(\boldsymbol{c}_{p} \lor x_{j})\Big) \thickapprox 1.
\]

Now, let $\A \in \class{K}$. As $\A \in \QQQ(\boldsymbol{Q}_{p}^{\gamma})$, we know that $\A$ is a subdirect product of algebras that are relatively subdirectly irreducible in $\QQQ(\boldsymbol{Q}_{p}^{\gamma})$. Since $\QQQ(\boldsymbol{Q}_{p}^{\gamma}) = \III\SSS\PPP\PPU(\boldsymbol{Q}_{p}^{\gamma})$, these algebras belong to $\III\SSS\PPU(\boldsymbol{Q}_{p}^{\gamma})$ and, therefore, are chains. Thus, $\A$ is a subdirect product of chains $\{ \C_{i} : i \in I \}$ in $\QQQ(\boldsymbol{Q}_{p}^{\gamma})$. Furthermore, as $\A$ validates the equation in the above display, so do the various $\C_{i}$. In view of Proposition \ref{Prop:RG-chains} and the fact that $\C_{i} \in \QQQ(\boldsymbol{Q}_{p}^{\gamma})$, this implies that
\[
\{ \C_{i} : i \in I \} \subseteq \III\SSS\PPU(\{ \boldsymbol{Q}_{p}^{0}, \dots, \boldsymbol{Q}_{p}^{n} \}) \subseteq \III\SSS\PPU(\boldsymbol{Q}_{p}^{n}).
\]
Thus, $\A$ is a subdirect product of members of $\QQQ(\boldsymbol{Q}_{p}^{n})$. Since $\QQQ(\boldsymbol{Q}_{p}^{n} ) = \QQQ(\Fm_{\class{K}}(\omega))$, we conclude that $\A \in \QQQ(\Fm_{\class{K}}(\omega))$. Thus, $\class{K}$ is generated as a quasivariety by $\Fm_{\class{K}}(\omega)$, as desired.
\end{proof}

The next result presents bases for the admissible rules on all the axiomatic  extensions of $\rgolog$.

\begin{Theorem}\label{Thm:SCbase-godel}
The following holds for every $r\in(0,1]$, $p\in[0,1)\cap\mathbb{Q}$, and $\gamma\in\omega+1$:
\benroman
\item A base for the admissible rules of\/ $\rgolog_{r}$ is given by the rules of the form $\boldsymbol{c}_{q} \lor z \rhd z$, for all $q \in [0, r) \cap \mathbb{Q}$;
\item A base for the admissible rules of \/ $\rgolog_{p}^{\gamma}$ is given by the rule $\boldsymbol{c}_{p} \lor z \rhd z$.
\eroman
\end{Theorem}

\begin{proof}
After Theorem \ref{Thm:RGA-SC},  the structural completion of  $\rgolog_{r}$ is algebrized by $\QQQ(\boldsymbol{Q}_{r})$ and   that of   $\rgolog_{p}^{\gamma}$ by $\QQQ(\boldsymbol{Q}_{p}^{\gamma})$. Thus, in order to obtain a base for admissible rules of $\rgolog_{r}$ and $\rgolog_{p}^{\gamma}$,  it suffices to find an axiomatization of $\QQQ(\boldsymbol{Q}_{r})$ and  $\QQQ(\boldsymbol{Q}_{p}^{\gamma})$   relative to $\VVV(\boldsymbol{Q}_{r})$ and $\VVV(\rgolog_{p}^{\gamma})$, respectively.   By Proposition \ref{Prop:RG-chains} and Theorem \ref{Thm:RGA-varieties}, the universal class of $\boldsymbol{Q}_{r}$ is axiomatized relative to $\VVV(\boldsymbol{Q}_{r})$ by $\boldsymbol{c}_{q}\not\thickapprox 1$ for all  $q \in [0, r) \cap \mathbb{Q}$. Using a similar argument as in  the proof of Theorem \ref{Thm:base-product},  $\QQQ(\boldsymbol{Q}_{r})$ is axiomatized relative to $\VVV(\boldsymbol{Q}_{r})$ by $\boldsymbol{c}_{q}\lor z \thickapprox 1 \Longrightarrow z \thickapprox 1$ for all  $q \in [0, r) \cap \mathbb{Q}$. Similarly, $\QQQ(\boldsymbol{Q}_{p}^{\gamma})$ is axiomatized relative to $\VVV(\boldsymbol{Q}_{p}^{\gamma})$ by $\boldsymbol{c}_{p}\lor z \thickapprox 1 \Longrightarrow z \thickapprox 1$.
\end{proof}

\section{Rational \L ukasiewicz logic}\label{Sec:9}

The lattice of axiomatic extensions of $\lulog$ is denumerable and its structure was completely described in \cite{Ko81}, see also \cite[Chpt.\ 8]{CiMuOt99}. On the other hand, the variety of MV-algebras is well known to be $\mathcal{Q}$-universal \cite{AdaDzi94}. We conclude this paper by showing that the addition of constants trivializes the lattice of extensions of $\rlulog$.

\begin{Theorem}
\label{thm_structural_completeness_RMV}
The logic $\rlulog$ has no proper consistent extensions and, therefore, is HSC.
\end{Theorem}

The rest of this section is dedicated to the proof of the theorem. 
Recall that $\rationalL$ denotes the subalgebra of $\standardL$ on the rational numbers in $[0,1]$
and ${\mathrmL}_{n+1}$ is the subalgebra on $\{0,\nicefrac{1}{n},\dots,\nicefrac{n-1}{n},n\}$.
The  algebra $\rationalL$ generates $\Alg{MV}$ as a quasivariety, 
as the following proposition shows:\footnote{In fact,
by \cite[Thm.~2.5]{Gp02}, $\III\SSS\PPU(\standardL)$ and $\III\SSS\PPU(\rationalL)$ coincide.}

\begin{Proposition}[\protect{\cite[Thm.\ 17]{Aguzzoli-Ciabattoni:Finiteness}}]
\label{prop_finiteness_in_inf_valued_Luk} 
Assume a MV-quasiequation $\Phi$
does not hold in $\Alg{MV}$ (equivalently, in $\standardL$). 
Then there is a natural number $n$ % that can be computed from the quasiequation,
such that $\Phi$ does not hold in ${\mathrmL}_{n+1}$.
\end{Proposition}

The details for the next paragraph can be found in \cite{Ha98}, see also \cite{Hanikova:ImplicitDefConstLuk}. 
Let $\alg{A}$ be an MV-algebra, 
$a$ an element of its universe, $\varphi(x_1,\dots,x_n)$ a term in the language of MV-algebras, and $1\leq i\leq n$. 
The equation $\varphi\thickapprox 1$ \emph{implicitly defines} $a$ in variable $x_i$ in $\alg{A}$  
provided that 
there are $c_1, \dots, c_n \in A$ such that $\varphi^\A(c_1, \dots, c_n) = 1^\alg{A}$ and, moreover, for every $c_1, \dots, c_n \in A$, 
\[
\text{if }v(c_1, \dots, c_n)=1^\alg{A}\text{, then }c_i =a.
\] 
The element $a$ is said to be implicitly definable in $\alg{A}$ if there is an equation that 
implicitly defines it.

\begin{Proposition}[\protect{\cite[Lems.\ 3.3.11 \& 3.3.13]{Ha98}}]
\label{prop_implicit_definability}
{\ }
\benroman
\item \label{item:imp-def-1} Each rational number $r\in (0,1)$ is implicitly definable in $\standardL$.
\item \label{item:imp-def-2} If $T\cup\{\alpha \thickapprox \beta\}$ is a finite set of equations in variables $\vec{x}$ and 
constants $\boldsymbol{c}_{r_1},\dots,\boldsymbol{c}_{r_k}$ in the language of rational MV-algebras, 
$T_0$ is a (finite) set of equations in the language of MV-algebras implicitly defining each $r_i$, $1\leq i\leq k$, 
in distinct variables $z_i$ (not among $\vec{x}$) in $\standardL$, 
and $T^\star, \alpha^\star,\beta^\star$ result from $T,\alpha,\beta$ by replacing the constant
$\boldsymbol{c}_{r_i}$ with variable $z_i$ respectively for $1\leq i\leq k$, then 
$T\models_{\standardLQ} \alpha\thickapprox \beta$ if and only if 
$T_0\cup T^\star \models_{\standardL}\alpha^\star\thickapprox \beta^\star$;
\item \label{item:imp-def-3} Under the conditions and notation of (ii), also 
$T\models_{\rationalLQ} \alpha\thickapprox \beta$ if and only if 
$T_0\cup T^\star \models_{\rationalL}\alpha^\star\thickapprox \beta^\star$.
\eroman
\end{Proposition} 

\begin{proof} Part (\ref{item:imp-def-1}) is \cite[Lem.\ 3.3.11]{Ha98} and part (\ref{item:imp-def-2})
is the proof method of \cite[Lem.\ 3.3.13]{Ha98}: the variables $z_i$, 
under the theory implicitly defining the finitely many rational constants, 
act semantically as the respective constants. 
Since the interpretation of constants is with the rationals, the argument 
can be carried out for the algebras $\rationalLQ$ and $\rationalL$, which justifies (\ref{item:imp-def-3}), 
even though not explicit in \cite{Ha98}.
\end{proof}

Notice that, in contrast to the situation in G\"odel or product algebras, the
interpretation of the constant $\boldsymbol{c}_r$ with the rational $r$ is the only possible
in the algebra $\standardL$; see under Thm.~26 in \cite{EsGiGoNo07} for an explicit mention of this fact.

\begin{Proposition}[\protect{\cite[Lem.\ 3.3.14]{Ha98}}, \protect{\cite[Prop.~24]{EsGiGoNo07}}]
\label{RMV_generated_real}
$\class{RMV} =\QQQ(\standardLQ)$.
\end{Proposition}

\begin{Lemma}
\label{RMV_generated_rational}
$\standardLQ$ and $\rationalLQ$ generate the same quasivarieties.
\end{Lemma}

% say that both have the canonical interpretation of rational constants. 
% Either this will be part of Section 3 or we need to say it here.

\begin{proof}
Let $\Phi = (\varphi_{1} \thickapprox \psi_{1} \boland \dots \boland \varphi_{n} \thickapprox \psi_{n}) \boto \varphi \thickapprox \psi$ be a quasiequation in the language of $\class{RMV}$. Moreover let 
$\boldsymbol{c}_{r_1},\dots,\boldsymbol{c}_{r_k}$ be the constants in $\Phi$,
and assume $T_0$ is a finite conjunction of equations in the language of MV-algebras that implicitly define $r_1,\dots,r_k$ in $\standardL$
in some pairwise distinct 
variables  $z_1,\dots,z_k$ not occurring in $\Phi$. (We further assume that any auxiliary variables in 
the defining equations for $r_i, r_j$, where $i\not= j$, are also distinct.)
Notice that $T_0$ exists by Proposition \ref{prop_implicit_definability}(\ref{item:imp-def-1}). Finally, let $\Phi^\star$ result from $\Phi$ by replacing the rational 
constants $\boldsymbol{c}_{r_1},\dots,\boldsymbol{c}_{r_k}$ with $z_1,\dots,z_k$. 
Notice that $\Phi^\star$ is an quasiequation in the language of MV-algebras.
% what exactly do we do about 0 and 1?

It is enough to show that $\rationalLQ\models\Phi$ implies 
$\standardLQ\models\Phi$. 
To that end, assume $\standardLQ\nvDash\Phi$; by Proposition \ref{prop_implicit_definability}(ii) 
this is equivalent to $\standardL\nvDash T_0 \boto \Phi^\star$. Moreover, using Proposition \ref{prop_finiteness_in_inf_valued_Luk}, there is a natural number $n$ 
 such that ${\mathrmL}_{n+1}\nvDash T_0 \boto \Phi^\star$.
Since each finite MV-chain is isomorphic to a subalgebra of $\rationalL$, we may conclude
that $\rationalL\nvDash T_0 \boto \Phi^\star$.
Finally, the latter is the case if and only if $\rationalLQ\nvDash\Phi$, applying 
Proposition \ref{prop_implicit_definability}(\ref{item:imp-def-3}). 
\end{proof}

\begin{Proposition}
\label{lm_diagram_lemma}
 Let $\alg{A}\in\Alg{RMV}$ be nontrivial. Then 
$\rationalLQ\leq\alg{A}$.
\end{Proposition}
\begin{proof}
Any nontrivial $\alg{A}\in\Alg{RMV}$ has a subalgebra that is a homomorphic image of 
$\rationalLQ$, and since the latter is simple and $\A$ lacks trivial subalgebras, 
we get $\rationalLQ\leq\alg{A}$.
\end{proof}

\begin{proof}[Proof of Theorem \ref{thm_structural_completeness_RMV}.]
We will show that \Alg{RMV} is a minimal quasivariety. 
By Proposition \ref{lm_diagram_lemma}, $\rationalLQ\leq\alg{A}$ 
whenever $\alg{A}\in\Alg{RMV}$ is nontrivial.
Thus, for any nontrivial class $\Alg{K}\subseteq\Alg{RMV}$, we have $\QQQ(\rationalLQ)\subseteq \QQQ(\Alg{K})$.
Combining Proposition \ref{RMV_generated_real} and Lemma \ref{RMV_generated_rational},
we get $\Alg{RMV}=\QQQ(\rationalLQ)$; since we know the latter to be included in $\QQQ(\Alg{K})$,
and given that $\Alg{K}$ was arbitrarily chosen, the minimality of $\Alg{RMV}$ follows.
\end{proof}

\paragraph{\bfseries Acknowledgements.}
The first and the third authors were supported by the research grant $2017$ SGR $95$ of the 
AGAUR from the Generalitat de Catalunya and  by the I+D+i research project PID2019-110843GA-I00 
\textit{La geometria de las logicas no-clasicas} funded by the Ministry of Science and Innovation of Spain. 
The second author was supported partly by the grant GA18-00113S of the 
Czech Science Foundation and partly by the long-term strategic
development financing of the Institute of Computer Science RVO:67985807.
 The third author was also supported by the \textit{Beatriz Galindo} grant BEAGAL\-$18$/$00040$ funded by the Ministry of Science and Innovation of Spain and by the grant CZ.$02$.$2$.$69$/$0$.$0$/$0$.$0$/$17$\_$050$/$0008361$, OPVVV M\v{S}MT, MSCA-IF \textit{Lidsk\'{e} zdroje v teoretick\'{e} informatice} funded by the Ministry of Education, Youth and Sports of the Czech Republic.

\bibliographystyle{plain}

\end{document}